\documentclass[aop,seceqn,citesort,MSNbibl,dvips]{arximspdf}

\doi{10.1214/10-AOP555}
\volume{39}
\issue{2}
\pubyear{2011}
\firstpage{439}
\lastpage{470}

\makeatletter

\newcommand{\eqref}[1]{(\ref{#1})}

\newtheorem{theorem}{Theorem}[section]
\newproclaim{rmk}{Remark}[section]
\newtheorem{lemma}{Lemma}[section]

\newcommand{\N}{\mathrm{N}}
\newcommand{\PP}{\mathrm{P}}
\newcommand{\Q}{\mathrm{Q}}
\newcommand{\E}{\mathrm{E}}
\newcommand{\F}{\mathrm{F}}

\makeatother

\begin{document}
\begin{frontmatter}

\title{An extension of the L{\'e}vy characterization to fractional
Brownian motion}
\runtitle{Characterization of fBm}

\begin{aug}
\author[A]{\fnms{Yuliya} \snm{Mishura}\thanksref{t1}\ead[label=e1]{myus@univ.kiev.ua}} and
\author[B]{\fnms{Esko} \snm{Valkeila}\corref{}\thanksref{t2}\ead[label=e2]{esko.valkeila@tkk.fi}}
\runauthor{Y. Mishura and E. Valkeila}
\affiliation{Kiev University and Aalto University}
\address[A]{Department of Mathematics\\
Kiev University\\
Volodomirska Street 64\\
01033 Kiev\\
Ukraine\\
\printead{e1}}
\address[B]{Department of Mathematics\\
\quad and Systems Analysis\\
Aalto University\\
P.O. Box 11100\\
FI-00076 Aalto\\
Finland\\
\printead{e2}}
\end{aug}

\thankstext{t1}{Supported in part by by the Suomalainen Tiedeakatemia.}
\thankstext{t2}{Supported by the Academy of Finland Grants 210465 and
212875.}

% HISTORY:
\received{\smonth{3} \syear{2008}}
\revised{\smonth{2} \syear{2010}}

% ABSTRACT
%
\begin{abstract}
Assume that $X$ is a continuous square integrable process with zero mean,
defined on some probability space $(\Omega, \F, \PP) $.
The classical characterization due to P. L{\'e}vy says that
$X$ is a Brownian motion if and only if $X$ and $X_t^2 - t $, $t\ge0,
$ are
martingales with respect to the intrinsic filtration $\F^ X$.
We extend this result to fractional Brownian motion.
\end{abstract}

% KEYWORDS
%
\begin{keyword}[class=AMS]
\kwd[Primary ]{60G15}
\kwd[; secondary ]{60E05}
\kwd{60H99}.
\end{keyword}
\begin{keyword}
\kwd{Fractional Brownian motion}
\kwd{L{\'e}vy theorem}.
\end{keyword}

\end{frontmatter}

%s1 ###
\section{Introduction}

In classical stochastic analysis, L{\'e}vy's characterization result for
standard Brownian motion is a fundamental result. We extend L{\'e}vy's
characterization result to fractional Brownian motion, giving three
necessary and sufficient properties for the process $X$ to be a fractional
Brownian motion. Fractional Brownian motion is a self-similar Gaussian
process with stationary increments. However, these two properties are not
explicitly present in the three conditions we shall give.

Fractional Brownian motion is a popular model in applied probability, in
particular, in teletraffic modeling and, to some extent, in finance.
Fractional Brownian motion is not a semimartingale and there has been much
research on how to define stochastic integrals with respect to fractional
Brownian motion. A large part of the developed theory depends on the
fact that
fractional Brownian motion is a Gaussian process. Since we want to prove
that $X$ is a special Gaussian process, we cannot use this machinery
for our
proof. L{\'e}vy's characterization result is based on It{\^o} calculus. We
cannot perform computations using the process $X$. Instead, we use the
representation of the process $X$ with respect to a certain martingale. In
this way, we can perform computations using methods from classical stochastic
analysis.

\subsection*{Notation and definitions}

We use the following notation: $\stackrel{L^p(\PP)}{\rightarrow}$
means convergence in the space $L^p(\PP)$, $\stackrel{\PP
}{\rightarrow}$
(resp., $\stackrel{\mathrm{a.s.}}{\rightarrow}$) means convergence in
probability
(resp., almost sure convergence) and $B(a,b)$ is the beta integral
$B(a,b) =
\int_0^1 x^{a-1}(1-x)^{b-1}\,dx $, defined for $a,b > 0 $. The notation $X_n
\le Y + o_\PP(1) $ means that we can find random variables $\epsilon_n
$ such
that $\epsilon_n = o_\PP(1) $ and $X_n \le Y + \epsilon_n $. If, in addition,
we have $X = \PP-\lim X_n $ in such a situation, then $X\le Y $.

If $M$ is a continuous square integrable martingale, then
the bracket of $M$ is denoted by $[M]$. Recall that, in this case, we have
\[
[M]_t = \PP-\lim_{|\pi^n| \to0 } \sum_{k=1} ^n
(M_{t_k}-M_{t_{k-1}})^2.
\]

\subsection*{Fractional Brownian motion}

A continuous square integrable centered process $X=(X_{t})_{t\geq0}$
with $%
X_{0}=0$ is a \textit{fractional Brownian motion} with
self-similarity index $H\in(0,1)$ if it is a Gaussian process with zero
mean and covariance function
%
%e1.1 ###
\begin{equation} \label{eq:fBm}
\E( X_{s}X_{t} ) =\tfrac{1}{2} (
t^{2H}+s^{2H}-|t-s|^{2H} ),\qquad s,t\geq0.
\end{equation}
If $X$ is a continuous Gaussian process with covariance
(\ref{eq:fBm}), then,
obviously, $X$~has stationary increments and $X$ is self-similar with
index $H
$. Mandelbrot named the Gaussian process $X$ from \eqref{eq:fBm}
\textit{fractional Brownian motion} and proved an important representation result
for fractional Brownian motion in terms of standard Brownian motion in
\cite{mvn}. For results concerning fractional Brownian motion before
Mandelbrot, we refer to \cite{mol}.

\subsection*{Characterization of fractional Brownian motion}

Throughout this paper, we work with special partitions. For $
t>0$, we put $t_k:= t\frac{k}{n} $, $k=0,\ldots, n$. Further, let $\F
^X$ be
the filtration generated by the process $X$.
Fix $H\in(0,1)$. Fractional Brownian motion has the
following three properties:

\begin{enumerate}[(a)]
\item[(a)] the sample paths of the process $X$ are $\beta$-H
{\"o}lder continuous for any $\beta\in(0,H)$;

\item[(b)] for $t>0$, we have
%
%e1.2 ###
\begin{equation} \label{eq:modified}
n^{2H-1}\sum_{k=1}^n (X_{t_k} - X_{t_{k-1}}) ^2 \stackrel{L ^1 (\PP
)}{%
\longrightarrow} t^{2H}
\end{equation}
as $n\to\infty$;

\item[(c)] the process
%
%e1.3 ###
\begin{equation} \label{eq:mol}
M_t = \int_0^ t s^{1/2 -H } (t-s) ^{1/2 - H} \,d X_s
\end{equation}
is a martingale with respect to the filtration $\F^X$.
\end{enumerate}

If the process $X$ satisfies (a), then we say that it is
\textit{H{\"o}lder up to} $H$. The property (b) characterizes
the \textit{weighted quadratic variation} of the process $X$ and the process $M$ in
(c) is the \textit{fundamental martingale} of $X$. It is a
martingale with the bracket $c_H t^{2-2H}$ for some constant $c_H$ and
is actually a time-changed Brownian motion, up to a constant. It
follows from
property (a) that the integral (\ref{eq:mol}) can be
understood as a Riemann--Stieltjes integral (see \cite{nvv} and Section
\ref{ip} for more details).

Fractional Brownian motion satisfies property (a):
from (\ref{eq:fBm}), we have that
\[
\E(X_t -X_s ) ^ 2 = (t-s)^{2H}.
\]
Since the process $X$ is a Gaussian process, we obtain from Kolmogorov's
theorem~\cite{ry}, Theorem I.2.1, page 26, that the process $X$ is
$\beta
$-H{\"o}lder
continuous with $\beta< H$. Fractional Brownian motion also satisfies
property (b). The proof of this fact is based on the
self-similarity and the ergodicity of the fractional Gaussian noise
sequence $Z_k:= X_k -X_{k-1}$, $k\ge1$. The fact that property (c)
holds for fractional Brownian motion was established in Molchan
\cite{mol} and
recently rediscovered by several authors (see \cite{nvv}).

We now
summarize our main result.
\begin{theorem}
\label{t:levy} Assume that $X$ is a continuous square
integrable centered process with $X_0=0$. Then, the following
properties are equivalent:

\begin{itemize}
\item the process $X$ is a fractional Brownian motion with
self-similarity index $H \in(0,1)$;

\item the process $X$ has properties \textup{(a)},
\textup{(b)} and \textup{(c)} for some $H\in(0,1)$.
\end{itemize}
\end{theorem}
\begin{rmk}
Theorem \ref{t:levy} appears in \cite{mis} with a
different proof.
\end{rmk}

\subsection*{Discussion}

If $H =\frac12 $, then assumption (c) means
that the process $X$ is a martingale. If $X$ is a martingale, then
condition (b) means that $X^2_t - t$ is a martingale. Hence, we
obtain the classical L{\'e}vy characterization theorem when $H=\frac12$.
Note that, in this case, property (a) follows from the fact
that $%
X$ is a standard Brownian motion.

Fractional Brownian motion $X$ also has the following property (see,
e.g., \cite{rogers}):
for $t > 0, $
%
%e1.4 ###
\begin{equation} \label{eq:1-H}
\sum_{k=1} ^n | X_{t_k} - X_ {t_{k-1}}| ^{{1/H}}
\stackrel{L^1 (\PP) }{\rightarrow} \E|X_1| ^{{1/H}} t
\end{equation}
as $n\to\infty$. To check that (\ref{eq:1-H}) holds for fractional
Brownian motion, similarly to~\eqref{eq:modified}, one can use
self-similarity and ergodicity of the fractional
Gaussian noise sequence. This provides another possibility to
generalize the
quadratic variation property of standard Brownian motion. However, it is
difficult to replace condition (b) by the condition (\ref%
{eq:1-H}).
\begin{rmk}
In the recent work of Hu et al. \cite{hns}, condition
(b) is replaced by the condition (\ref{eq:1-H}), with the
additional assumption that $[M]$ is absolutely continuous with respect to
the Lebesgue measure for $H > \frac12 $. The authors\vspace*{1pt} show that
conditions (a), (c) and (\ref{eq:1-H}) also characterize
fractional Brownian motion. In our work, we do not suppose the absolute
continuity of $[M]$, but prove it under other assumptions; however, we
restrict ourselves to (b).
\end{rmk}

In the next section\vspace*{2pt} we give one auxiliary result. The rest of the
paper is devoted to the proof of the main result, first for $H > \frac
12 $ and then for $%
H<\frac12 $.

%s2 ###
\section{Auxiliary result}

%s2.1 ###
\subsection{Martingales and random variables}

In the proof, we will use random variables which are final values
of martingales of a special type. All martingales vanish at zero.

Two continuous martingales $M,N$ are (\textit{strongly})
\textit{orthogonal} if $[M,N]=0$; we write this as $M\perp N$. Integration by parts
gives that for such $M,N$, the product $MN$ is a local martingale and
it then
has a bracket $[MN]$. We use the notation $N\cdot M$ for the stochastic
integral of $N$ with respect to $M\dvtx (N\cdot M)_t = \int_0^t
N_s\,dM_s$. Let
$M$ be a continuous martingale. Put $I_2(M) _t := (M\cdot M)_t = \int_0^t
M_s \,dM_s $.

Let $0<a<b<t$ and suppose that $p,q$ are deterministic continuous
functions. Define martingales $N$ and $\tilde{N}$ by $N_{s}=%
\int_{0}^{s}p_{u}1_{(0,a]}(u)\,dM_{u}$ and $\tilde{N}_{s}=%
\int_{0}^{s}q_{u}1_{(a,b]}(u)\,dM_{u}$, respectively. The martingales $N$
and $\tilde{N}$
are orthogonal by construction and hence their product\vspace*{1pt} is a martingale.
Note that $N_{s}\tilde{N}_{s}=0$ whenever $s\leq a$ and $N_{s}\tilde
{N}%
_{s}=N_{a}\tilde{N}_{s}$ for $s>a$. The bracket of the martingale
$N\tilde{N}
$ is $[N\tilde{N}]_{s}=0$ whenever $s\leq a$ and $[N\tilde
{N}]_{s}=N_{a}^{2}[%
\tilde{N}]_{s}$ for $s>a$.

For orthogonal martingales, we have following lemma,
which we will use in our proof.
\begin{lemma}
\label{le:img} Assume that $(M^{n,k}_t)_{t\geq0}$
is a double array of
continuous square integrable martingales with the properties:

\begin{longlist}[(iii)]
\item[(i)] for $n$ fixed and $k\neq l$, $M^{n,k}$
and $%
M^{n,l}$ are orthogonal martingales;

\item[(ii)] for any $t\geq0$, $\sum_{k=1} ^{k_n}
[M^{n,k}]_t \le C $, where $C
$ is a constant;

\item[(iii)] for any $t\geq0$, $\max_ k [M^{n,k}]_t
\stackrel{\PP}{\to} 0 $ as $
n\to\infty$,
\end{longlist}
where $1\le k\le k_n$ and $k_n\to\infty$ as $n\to
\infty$;
then, for any $t\geq0$,
%
%e2.1 ###
\begin{equation} \label{eq:img}
\sum_{k=1} ^{k_n} I_2(M^{n,k})_t \stackrel{L^2(\PP)}{\to} 0
\end{equation}
as $n \to\infty$.
\end{lemma}
\begin{pf}
Since the martingales $M^{n,k}$ are pairwise orthogonal,
when $n$ is fixed, the same is true for
the iterated integrals $I_2(M^{n,k})$. Recall
\cite{ck}, Theorem 1, page~354, which states that
$\E( I_2(M^{n,k})_t ) ^2 \le B_{2,2} \E[M^{n,k}]^2 _t$.
Here, $B_{2,2}$ is constant independent of $n,t$ and $k$, and
this, together with property (ii), gives that the iterated
integrals $I_2(M^{n,k})$ are square
integrable. Hence, by the orthogonality of the iterated integrals, we have
\[
\E\Biggl( \sum_{k=1}^{k_n} I_2(M^{n,k})_t \Biggr) ^2 = \sum
_{k=1}^{k_n}\E( I_2(M^{n,k})_t ) ^{2}.
\]
However,
\[
\sum_{k=1}^{k_n}[M^{n,k}]^2 _t \le\max_k [M^{n,k}]_t \sum
_{k=1}^{k_n}[M^{n,k}] _t \stackrel{\PP}{\to} 0
\]
as $n\to\infty$. The claim \eqref{eq:img} now follows since $\max_k
[M^{n,k}]_t \le\sum_{k=1}^{k_n}[M^{n,k}] _t$ and this, together with
property (ii), gives
\[
\sum_{k=1}^{k_n}[M^{n,k}]^2 _t \le\max_k [M^{n,k}]_t \sum
_{k=1}^{k_n}[M^{n,k}] _t\le C^2.
\]
\upqed\end{pf}

%s2.2 ###
\subsection{A consequence of \textup{(b)}}
\label{ss:f-up-b}

We now fix $t$ and let $\mathcal{R }_t := \{ s \in[0,t] \dvtx
\frac{s}{t} \in\Q\} $. Note that the set $\mathcal{R}_t$ is dense on
the interval $[0,t]$. Now, also fix $s\in\mathcal{R}_t $ and let
$\tilde n =
\tilde n(s) $ be a subsequence of $n\in\N$ such that $\tilde n
\frac{s}{t} \in\N$. Put $\Delta X_{t_{k,n}} := X_{t_k
} - X_{t_{k-1}} $.

The next lemma opens the way to bound from below and above the bracket
$[M]$ on $[0, T]$ for any $T>0$ and this goal will be achieved in
Section \ref{ss:end-of-proof}.
\begin{lemma}
\label{le:subqv} Fix $t> 0$, $s\in\mathcal{R }_ t$ and
suppose that
$\tilde n \frac{s}{t}\in\N$ and $\tilde n
\to\infty$. Then,
\[
\tilde n ^{2H-1} \sum_{k=\tilde n {s/t} +1} ^{\tilde n} (
\Delta X_{t_{k,\tilde n}} ) ^2 \stackrel{L^1(\PP)}{%
\longrightarrow} t ^{2H-1}(t-s).
\]
\end{lemma}
\begin{pf}
We have that
\begin{eqnarray*}
&&\tilde{n}^{2H-1}
\sum_{k=1}^{\tilde{n}{s/t}} (\Delta
X_{t_{k,\tilde n}} )^{2} \\
&&\qquad =  \tilde{n}^{2H-1}\sum_{k=1}^{\tilde
{n}{s/t}} (\Delta
{X_{s_{k,\tilde{n} {s/t}}}})^{2}\\
&&\qquad =  \biggl(\frac{t}{s} \biggr)^{2H-1}\cdot\biggl({\tilde{n} \frac
{s}{t}} \biggr)^{2H-1}
\sum_{k=1}^{\tilde{n}{s/t}} (\Delta
{X_{s_{k,\tilde{n} {s/t}}}})^{2}
\stackrel{L^1(\PP)}{\longrightarrow}
s^{2H}\cdot\biggl(\frac{t}{s} \biggr)^{2H-1} \\
&&\qquad = st^{2H-1}.
\end{eqnarray*}
Since $\tilde{n}^{2H-1} \sum_{k=1}^{\tilde{n}} (\Delta
X_{t_{k,\tilde n}} )^{2}\stackrel{L^1(\PP)}{\longrightarrow}
t^{2H}$, we obtain
the proof.
\end{pf}

In what follows, we shall write $n$ for $\tilde n$ and $t_k$
for $t\frac{k}{n}$.

%s2.3 ###
\subsection{Some representation results}
\label{ip}

We shall use the following notation. Let
$Y_t = \int_0^t s^{1/2 -H} \,dX_s$.
We then we have $X_t = \int_0^t s^{H-1/2}\, dY_s $ and can write the
fundamental martingale $M$ as
%
%e2.2 ###
\begin{equation} \label{eq:m}
M_t = \int_0^t (t-s) ^{1/2 -H} \,dY_s.
\end{equation}

We also work\vspace*{1pt} with the martingale $W_t = \int_0^t
s^{H-1/2}\,dM_s$. We have $[W]_t = \int_0^t s^{2H-1}\,d[M]_s $ and
$[M]_t =
\int_0^ts^{1-2H}\, d[W]_s$.

The equation (\ref{eq:m}) is a generalized Abel integral
equation and the process $Y$ can be expressed in terms of the process $M$:
%
%e2.3 ###
\begin{equation} \label{eq:aym}
Y_t = \frac{1}{\Gamma(H+1/2 )\Gamma(3/2 -H)} \int_0^t (t-s)^{H-
1/2 } \,dM_s.
\end{equation}

Note that all of the integrals can be understood as pathwise
Riemann--Stieltjes integrals (see \cite{nvv}).

%s3 ###
\section{\texorpdfstring{Proof of Theorem \protect\ref{t:levy}: $H >
\frac12$}{Proof of Theorem 1.1: $H > \frac12$}}

%s3.1 ###
\subsection{Basic representation}

We shall now prove that $M$ is a martingale with a bracket
$c_Ht^{2-2H}$ for some constant $c_H$ and this, together with Lemma \ref{le:repXW},
will give that $X$ is a fractional Brownian motion with index $H$.

We shall use the following modified representation result
between $X$ and~$M$.
\begin{lemma}
\label{le:repXW} Assume that $H> \frac12 $ and that
properties \textup{(a)}
and \textup{(c)} hold. Then, the process $X$ has the representation
%
%e3.1 ###
\begin{equation} \label{eq:repXM}
X_t = \frac{1}{B_1} \int_0 ^t \biggl(\int_u^t s ^{H -1/2 } (
s-u )^{H-3/2}\,ds \biggr) \,d M_u
\end{equation}
with $B_1 = B (H-\frac12 , \frac32 -H ) $.
\end{lemma}
\begin{pf}
Integration by parts in
(\ref{eq:aym}) gives
\[
Y_t = \frac{1}{B_1} \int_0^t (t-s)^{H-3/2}M_s\,ds.
\]
Next, by using integration by parts and Fubini's
theorem, we obtain
\begin{eqnarray*}\quad
X_t &=& \int_0^t s^{H-{1/2}}\,dY_s \\
&=&
t^{H-{1/2}}Y_{t}- \biggl(H-\frac{1}{2} \biggr)
\int_{0}^{t}s^{H-{3/2}}Y_{s}\,ds \\
& = & \frac{t^{H-{1/2}}}{B_{1}}\int_{0}^{t}(t-s)^{H-
{3/2}}M_{s}\,ds\\
&&{} -
\frac{H-{1/2}}{B_{1}} \int_{0}^{t}s^{H-{3/2}}
\int_{0}^{s}(s-u)^{H-{3/2}}M_{u}\,du\,ds \\
& = & \frac{t^{H-{1/2}}}{(H-{1/2})B_{1}}\int
_{0}^{t}(t-s)^{H-{1/2}}\,dM_{s}\\
&&{}-
\frac{1}{B_1}\int_0^t
s^{H-{3/2}}\int_{0}^{s}(s-u)^{H-{1/2}}\,dM_{u}\,ds\\
& = & \frac{t^{H-{1/2}}}{(H-{1/2})B_{1}}\int
_{0}^{t}(t-s)^{H-{1/2}}\,dM_{s}\\
&&{} -
\frac{1}{B_{1}}\int_{0}^{t} \biggl[\int_{u}^{t}s^{H-
{3/2}}(s-u)^{H-{1/2}}\,ds \biggr]\,dM_{u}\\
& = & \frac{1}{B_{1}}\int_{0}^{t} \biggl[\frac{t^{H-{1/2}}}{H-
{1/2}}(t-u)^{H-{1/2}}\\
&&\hspace*{33.6pt}{}-
\int_{u}^{t}s^{H-{3/2}}(s-u)^{H-{1/2}}\,ds \biggr]\,dM_{u}\\
& = & \frac{1}{B_{1}}\int_{0}^{t} \biggl[\int_{u}^{t}s^{H-
{1/2}}(s-u)^{H-{3/2}}\,ds \biggr]\,dM_u.
\end{eqnarray*}
This proves claim (\ref{eq:repXM}).
\end{pf}

Our plan is now as follows: we will attempt to prove that $M$ is a
martingale with the bracket $C_Ht^{2-2H}$
and this, together with Lemma \ref{le:repXW}, will give that $X$ is a
fractional Brownian motion with parameter $H$.

%s3.2 ###
\subsection{The basic estimation}
\label{ss:basic}

We can assume that the processes $M$, $W$, $
[M]$ and $[W]$ are bounded with a deterministic constant $L$. If this
is not
the case, then consider a stopping time $\tau$,
\[
\tau= \inf\{ s\dvtx |M_s| \ge L \mbox{ or } |W_s| \ge L \mbox{ or }
\lbrack
M]_s \ge L \mbox{ or } \lbrack W]_s\ge L \}.
\]
Note that $\tau$ is independent of the partition $(t^n_k)$, $k=0,
\ldots, n$%
, and hence we have
\[
1_{\{\tau\ge t \} } n^{2H-1}\sum_{k=1}^n (\Delta
X_{t_{k, n}} )^{2} \stackrel%
{\PP}{\longrightarrow} 1_{\{\tau\ge t \} }t^{2H}.
\]
Given $\epsilon> 0, $ take $L$ big enough such that $\PP( \tau<
t )
< \epsilon$. Since our asymptotic results concern
convergence in probability, it is enough to prove them only in the set
$\{ \tau\ge t\}$. We do not write the stopping
time $\tau$ or the indicator $1_{\{\tau\ge t \} }$ explicitly in the proof
below.

We want to use the expression
\[
n ^{2H-1} \sum_{k=n {s/t} +1}^{ n} (\Delta
X_{t_{k, n}} )^{2}
\]
to obtain estimates for the increment of the bracket $[M]$, with the
help of
(\ref{eq:repXM}).

Use (\ref{eq:repXM}) to obtain
%
%e3.2 ###
\begin{equation}
\label{eq:rmk}
\Delta
X_{t_{k, n}} = \frac{1}{B_1} \biggl( \int_0 ^{t_{k-1}
}f^t_k (s ) \,dM_s + \int_{t_{k-1}} ^{t_k} g^t_k (s )\,dM_s
\biggr),
\end{equation}
where we have used the notation
%
%e3.3 ###
\begin{equation} \label{eq:f}
f_k^t (s) := \int_{t_{k-1}}^{t_k} u^{H-1/2} (u-s)^{H-3/2}\,du
\end{equation}
and
\[
g^t_k(s) := \int_s^{t_k}u^{H-1/2}(u-s)^{H-3/2}\,du.
\]
Rewrite the increment of $X$ as
%
%e3.4 ###
\begin{eqnarray} \label{eq:arr}\qquad
\Delta
X_{t_{k, n}} &=& \frac{1}{B_1} (I_k^{n,1} + I_k^{n,2} +
I_k^{n,3} ) \nonumber\\[-8pt]\\[-8pt]
:\!&=&
\frac{1}{B_1} \biggl( \int_0 ^{t_{k-2}} f^t _k (s)\,dM_s + \int_{t_{k-2}
}^{t_{k-1}} f^t_k(s)\,dM_s + \int_{t_{k-1}}^{t_k} g^t _k (s)\,dM_s \biggr).
\nonumber
\end{eqnarray}
We need such a decomposition because the behavior of the kernels in
the integrands is different
for different arguments. Now, we intend to use this decomposition and
to show that the sequence $n ^{2H-1} \sum_{k=n {s/t} +1}^{ n}
(\Delta
X_{t_{k, n}} )^{2}$ verifies relation (e) from Section
\ref{ss:end-of-proof}.
In order to do this, we use Lemma \ref{le:repXW}, decompose the
increment $\Delta X_{t_{k, n}}$
according to \eqref{eq:arr} into several terms and apply It{\^o's}
formula to the square of the increments. We then try to find
asymptotically nontrivial terms and terms of order $o_\PP(1)$, and
nontrivial terms must be of the form that will be appropriate for
finding the bounds for $[M]$. Even at this point, we can note that the
nontrivial terms will appear when we consider sums of the form $n
^{2H-1} \sum_{k = n {s/t} + 2 } ^{n} \int_0 ^{t_{k-2}}
(f^t_k(u))^2 \,d[M]_u$, etc. So, at first, we estimate the sums with such
a form and only then consider the remainder terms.

We note that the random variables $I^{n,j}_k$ are the final values at
moment $t$ of the
martingales $\int_0^{t_{k-2}\wedge v }
f^t_k(u)\,dM_u $, $ \int_{t_{k-2}\wedge v}^{t_{k-1}\wedge v}
f^t_k(u)\,dM_u $ and $ \int_{t_{k-1}\wedge v}^{t_k\wedge v}
g^t_k(u)\,dM_u $, $0\leq v\leq t$, respectively. By construction,
these martingales are orthogonal.

Next, the following upper bound holds for the functions $
f^t_k$:
%
%e3.5 ###
\begin{equation} \label{eq:f-up}
f^t_k (s)
% & \int_{t_{k-1}}^{t_k} u^{H-1/2} (u-s)^{H-3/2}\,du
\le t_k^{H-1/2} ( t_{k-1}-s )
^{H-3/2}\frac{t}{n};
\end{equation}
note that this estimate is finite (not bounded) for $s\in[0, t_{k-1} )
$ and bounded for
$s\in[0,t_{k-2}]$.
Further, we need the following technical result.
\begin{lemma}\label{le:app-1}
For $u< s$, we have
%
%e3.6 ###
\begin{eqnarray}\label{eq:app-1-a}
\sum_{k=n{s/t}+2}^{n}(t_{k-1}-u)^{2H-3}
&\le& \biggl(s+\frac{t}{n}
-u\biggr)^{2H-3}\nonumber\\[-8pt]\\[-8pt]
&&{} + \frac{n}{(2-2H)t}
\biggl(s+\frac{t}{n} -u\biggr)^{2H-2}\nonumber
\end{eqnarray}
and for $u\le t_i$, we have
%
%e3.7 ###
\begin{equation}\label{eq:app-1-b}\quad
\sum_{k=i+2}^n (t_{k-1}-u)^{2H-3} \le(t_{i+1}-u)^{2H-3} + \frac
{n}{(2-2H)t}(t_{i+1}-u)^{2H-2}.
\end{equation}
\end{lemma}
\begin{pf}
For $u < s$, we have
\begin{eqnarray*}
&&
\sum_{k=n{s/t}+2}^{n}(t_{k-1}-u)^{2H-3}\\
&&\qquad = \biggl( s + \frac
{t}{n} -u \biggr) ^{2H-3} + \frac{n}{t}\sum_{k=n
{s/t}+3}^{n}(t_{k-1}-u)^{2H-3}\frac{t}{n}\\
&&\qquad \le \biggl( s + \frac{t}{n} -u \biggr) ^{2H-3} + \frac{n}{(2-2H)t}
\biggl(s+\frac{t}{n} -u\biggr)^{2H-2}
\end{eqnarray*}
by estimating the second sum in the first line from above by the
integral. This proves
(\ref{eq:app-1-a}). Inequality (\ref{eq:app-1-b}) is proved in the
same way.
\end{pf}

We can now give two-sided bounds for the brackets of the martingales in
\eqref{eq:arr}. As was mentioned before, these brackets give rise to
nontrivial terms in our estimates.
\begin{lemma}
\label{le:y1}
Fix $t> 0$ and $s\in\mathcal{R }_ t$,
and let $
\tilde n$ be such that $\tilde n \frac{s}{t}\in\N$ and $\tilde n
\to\infty$ (we write $n$ instead of $\tilde n$ in what follows).
Then, there
exist two constants, $C_1$, $C_2> 0$, such that
%
%e3.8 ###
\begin{eqnarray} \label{eq:y1}\qquad
C_1 t ^{2H-1}\int_{s-{t/n}}^{t-2{t/n}} u^{2H-1} \,d[M]_ u
&\le&
n ^{2H-1} \sum_{k = n {s/t} + 2 } ^{n} \int_0 ^{t_{k-2}}
(f^t_k(u))^2 \,d[M]_u \nonumber\\[-8pt]\\[-8pt]
& \le& C_2t ^{4H-2} ([M]_t - [M]_s ) + o_\PP(1).\nonumber
\end{eqnarray}
\end{lemma}
\begin{pf} We will not write the constants explicitly.

\textit{Upper bound in} (\ref{eq:y1}).
First, we estimate
\[
i^{n}:=n^{2H-1}\sum_{k=n{s/t}+2}^{n}\int_{0}^{t_{k-2}}
(f_k^t(u))^{2}\,d[M]_{u}
\]
from above. From (\ref{eq:f-up}), we obtain
the following estimate for $i^{n} $:
%
%e3.9 ###
\begin{equation}\label{eq:est-1}
i^{n}
\le
n^{2H-3}t^{2H+1}\sum_{k=n{s/t}+2}^{n}
\int_{0}^{t_{k-2}}(t_{k-1}-u)^{2H-3} \,d[M]_{u}.
\end{equation}

We can assume that $0<s<t$ and $2\leq n
\frac{s}{t}\leq n-4$, and rewrite the estimate in (\ref{eq:est-1}) as
%
%e3.10 ###
\begin{eqnarray} \label{eq:est-2}\qquad
{i}^{n}
& \le& n^{2H-3}t^{2H+1}\Biggl(\sum_{i=1}^{n{s/t}} \sum_{k=n
{s/t}+2}^{n}
+\sum_{i=n{s/t}+1}^{n-2} \sum_{k=i+2}^{n}\Biggr)\nonumber\\
&&{}\times
\int_{t_{i-1}}^{t_i} (t_{k-1}-u)^{2H-3}\,
d[M]_{u} \nonumber\\[-8pt]\\[-8pt]
& = & n^{2H-3}t^{2H+1}\sum_{i=1}^{n{s/t}}
\int_{t_{i-1}}^{t_i}
\Biggl(\sum_{k=n{s/t}+2}^{n}(t_{k-1}-u)^{2H-3}\Biggr)\,d[M]_{u} \nonumber\\
& &{} + n^{2H-3}t^{2H+1}\sum_{i=n{s/t}+1}^{n-2}\int_{t_{i-1}}^{t_i}
\Biggl(\sum_{k=i+2}^{n} (t_{k-1}-u)^{2H-3}\Biggr)\,d[M]_{u}.\nonumber
\end{eqnarray}
We estimate the first term in the last equation in (\ref{eq:est-2}) by
(\ref{eq:app-1-a}):
%
%e3.11 ###
\begin{eqnarray} \label{eq:est-2d}
R^t_n
:\!&= & n^{2H-3}t^{2H+1}\sum_{i=1}^{n{s/t}}\int
_{t_{i-1}}^{t_i}\sum_{k=n{s/t}+2}^{n}
(t_{k-1}-u)^{2H-3}\,d[M]_u \nonumber\\
& \le& t^{2H-1} \int_0^s \biggl(t^2(ns + t -nu)^{2H-3}\\
&&\hspace*{47pt}{} + \frac{t}{2-2H}(ns+t-u)^{2H-2} \biggr)\,d[M]_u.\nonumber
\end{eqnarray}
Note that $(ns+t-nu)^{2H-3}$ and $(ns+t-nu)^{2H-2}$ are bounded and
both converge to $0$ as $n\to\infty$.
So, $ R^t_n= o_\PP(1)$, by the dominated convergence theorem.

For the second term in the last equation of (\ref{eq:est-2}),
we obtain, from (\ref{eq:app-1-b}), using the estimate
$(t_{i+1}-u)^{H-1/2} \le(\frac{t}{n} ) ^{H-1/2} $
and summing,
\begin{eqnarray*}
&&n^{2H-3}t^{2H+1}\sum_{i={ns/t}+1}^{n}
\int_{t_{i-1}}^{t_{i}}\biggl[(t_{i+1}-u)^{2H-3}+
\frac{n}{(2-2H)t} (t_{i+1}-u)^{2H-2}\biggr]
\,d[M]_{u}\\
&&\qquad\le
c_Ht^{4H-2}([M]_{t}-[M]_{s}).
\end{eqnarray*}

Hence, we have proven the upper bound (\ref{eq:y1})
and have
\[
i^{n}\leq
c_Ht^{4H-2}([M] _{t}-[M] _{s}) + o_\PP(1).
\]

\textit{Lower bound in} (\ref{eq:y1}).
We complete the proof of Lemma \ref{le:y1} by giving the lower bound.
%Recall first that
%$$
%i^{n}:=n^{2H-1}\sum_{k=n{s/t}+2}^{n}\int_{0}^{t_{k-2}}
%(f_k^t(u))^{2}d[M]_{u},
%$$
From the definition of $i^{n}$, we easily obtain a lower estimate:
%
%e3.12 ###
\begin{equation}\label{eq:l-1}
i^{n}\ge n^{2H-1}\sum_{k=n{s/t}+2}^{n}\int_{t_{k-3}}^{t_{k-2}}
(f_k^t(u))^{2}\,d[M]_{u}.
\end{equation}

%Recall also that
%$f^t_k(u) = \int_{t_{k-1}} ^{t_k}v^{H-1/2}(v-u)^{H-3/2}\,dv$
Further, for $u\in(t_{k-3},t_{k-2})$, $v\in(t_{k-1},t_k)$, we
have $v-u\le\frac{3}{n}t$, $u < v$ and we get the estimate
%
%e3.13 ###
\begin{equation}\label{eq:f-down}
(f_k^t (u ) )^2 \geq3^{2H-3}t^{2H-1}n^{1-2H}u^{2H-1}.
\end{equation}

We use (\ref{eq:f-down}) in the lower bound \eqref{eq:l-1} to obtain
\begin{eqnarray*}
i^{n} &\ge& 3^{2H-3}t^{2H-1} \sum_{k=n{s/t}+2}^{n}\int
_{t_{k-3}}^{t_{k-2}}u^{2H-1}\,d[M]_u\\
&=& 3^{2H-3}t^{2H-1}\int_{s-{t/n}}^{t-2{t/n}} u^{2H-1 } \,d[M]_u
\end{eqnarray*}
and this gives the lower bound in (\ref{eq:y1}).
The proof of Lemma \ref{le:y1} is now complete.
%\rightqed
\end{pf}
\begin{rmk} Clearly, we can rewrite $i^{n}$ similarly to \eqref
{eq:est-2} as
%
%e3.14 ###
\begin{eqnarray} \label{eq:est-2a}\hspace*{28pt}
{i}^{n}
& = & n^{2H-3}t^{2H+1}\Biggl(\sum_{i=1}^{n{s/t}} \sum_{k=n{s/t}+2}^{n}
+\sum_{i=n{s/t}+1}^{n-2} \sum_{k=i+2}^{n}\Biggr)
\int_{t_{i-1}}^{t_i} (f_k^t(u))^{2}\,
d[M]_{u} \nonumber\\
& = & n^{2H-3}t^{2H+1}\sum_{i=1}^{n{s/t}}
\int_{t_{i-1}}^{t_i}
\sum_{k=n{s/t}+2}^{n}(f_k^t(u))^{2}\,d[M]_{u} \\
& &{} + n^{2H-3}t^{2H+1}\sum_{i=n{s/t}+1}^{n-2}\int_{t_{i-1}}^{t_i}
\sum_{k=i+2}^{n} (f_k^t(u))^{2}\,d[M]_{u} \nonumber
\end{eqnarray}
and obtain from \eqref{eq:f-up}, and similarly to \eqref{eq:est-2d}, that
%
%e3.15 ###
\begin{equation} \label{eq:est-2g}
n^{2H-3}t^{2H+1}\sum_{i=1}^{n{s/t}}
\int_{t_{i-1}}^{t_i}
\sum_{k=n{s/t}+2}^{n}(f_k^t(u))^{2}\,d[M]_{u}
\stackrel{\PP}{\rightarrow}0,
\end{equation}
change summation indices for further convenience and deduce from \eqref
{eq:est-2a}, \eqref{eq:est-2g} that
%
%e3.16 ###
\begin{eqnarray} \label{eq:est-2f}\quad
&&\PP-\lim_{n\rightarrow\infty}{i}^{n}\nonumber\\[-8pt]\\[-8pt]
&&\qquad= \PP-\lim_{n\rightarrow\infty} n^{2H-3}t^{2H+1}\sum_{k=n
{s/t}+1}^{n-2}\int_{t_{k-1}}^{t_k}
\sum_{i=k+2}^{n} (f_i^t(u))^{2}\,d[M]_{u}.\nonumber
\end{eqnarray}
\end{rmk}

We now return to \eqref{eq:arr}, take the bracket of the next term and
so estimate the term
\[
\int_{t_{k-2}} ^{t_{k-1}}(f^t_k (s))^2\,d[M]_s.
\]
\begin{lemma}
\label{le:y2} There exists a constant $C_{3}>0$ such that
%
%e3.17 ###
\begin{equation} \label{eq:f-y-12}\quad
n^{2H-1} \sum_{k=n{s/t}+2}^{n} \int
_{t_{k-2}}^{t_{k-1}}(f^t_k(u))^{2}\,d[M]_u\leq C_{3}t^{4H-2}([M]_{t}-[M]_{s}).
\end{equation}
\end{lemma}
\begin{pf}
We have the following upper estimate for the function $f^t_k$:
\begin{eqnarray*}
f^t_k (u) & \le& t_k^{H-1/2} \int_{t_{k-1}}^{t_k}
(v-u)^{H-3/2} \,dv \\
& = & \frac{1}{H-1/2} t_k^{H-1/2} \bigl( ( t_k-u
)^{H-1/2} - ( t_{k-1}-u ) ^{H-1/2} \bigr) \\
& \le& \frac{1}{H-1/2}t^{H-1/2} \biggl(\frac{t}{n} \biggr)
^{H-1/2}.
\end{eqnarray*}
This gives the claim \eqref{eq:f-y-12}.
\end{pf}

The last estimate for nontrivial terms in \eqref{eq:arr} concerns the
terms of the form
\[
\int_{t_{k-1}}^{t_k} (g^t_k(s))^2 \,d[M]_s.
\]
\begin{lemma}
\label{le:y3} There exists a constant $C_4$ such that
%
%e3.18 ###
\begin{equation} \label{eq:y3}\quad
n ^{2H-1} \sum_{k = n {s/t} + 1 }
^{n}\int_{t_{k-1}}^{t_k}(g^t_k(u))^2 \,d[M]_u \le C_4 t ^{4H-2} (
[M]_t -
[M]_s ).
\end{equation}
\end{lemma}
\begin{pf} We have that
\begin{eqnarray*}
g^t_k(z) &= &
\int_{z}^{t_{k}}v^{H-{1/2}}(v-z)^{H-{3/2}}\,dv\leq
(t_{k})^{H-{1/2}}
\frac{(t_{k}-z)^{H-{1/2}}}{H-{1/2}}\\
&\leq& C (t_{k})^{H-{1/2}} \biggl(\frac{t}{n} \biggr)^{H-
{1/2}}\leq
Ct^{2H-1} \biggl(\frac{1}{n} \biggr)^{H-{1/2}}.
\end{eqnarray*}
This gives the claim \eqref{eq:y3}.
\end{pf}
%

%s3.3 ###
\subsection{The $o_\PP(1)$ terms}
\label{ss:remainder}
%into several terms and apply It{\^o} formula to the square of the
%increments. Then we try to find asymptotically non-trivial terms and
%terms of order $o_\PP(1)$,

We shall now prove that after the decomposition of the increment
$\Delta X_{t_{k, n}}$
according to \eqref{eq:arr}, taking the square of this increment and
applying It{\^o's} formula to the decomposition, all the terms
except the three brackets of the martingales become asymptotically trivial.
In this order, we take the terms of the form $(I^{n,j}_k)^2 $, $
j=1,2,3$,
decompose them by It{\^o's} formula on the bracket and martingale part
and also prove that
the terms containing the cross products $I^{n,i}_kI^{n,j}_k $,
$i\neq j$, are asymptotically trivial.
More exactly, It{\^o's} formula implies that
\[
(I^{n,1}_k)^2 = \int_0^{t_{k-2}}(f^t_k(v))^2 \,d[M]_v + 2\int_0 ^{t_{k-2}}
f^t_k(u) \biggl( \int_0^u f^t_k(v) \,dM_v \biggr) \,dM_u.
\]
We shall show that
%
%e3.19 ###
\begin{equation} \label{eq:y-e-f-1}
n^{2H-1} \sum_{k=n{s/t} +2} ^n \int_0 ^{t_{k-2}} f^t_k(u) \biggl(
\int
_0^u f^t_k(v) \,dM_v \biggr) \,dM_u \stackrel{\PP}{\rightarrow} 0
\end{equation}
as $n\to\infty$. Clearly, it is sufficient to consider the sums of
the form
\[
S^{n} = n^{2H-1} \sum_{k=3}^{n}
\int_{0}^{t_{k-2}} \biggl(\int_{0}^{u} f_{k}^{t}(s) \,d M _{s} \biggr)
f_{k}^{t}(u)\, d M_{u},
\]
(note that $n\frac{s}{t}\geq1$) since the
sums
\[
\sum_{k=3}^{n{s/t}+1}
\int_{0}^{t_{k-2}} \biggl( \int_{0}^{u} f_{k}^{t}(s) \,d M _{s} \biggr)
f_{k}^{t}(u)\,d M_{u}
\]
for $n\frac{s}{t}\geq2$ can be considered in
a similar way. We rewrite $S^{n}$ as
\begin{eqnarray*}
S^{n} &=& n^{2H-1} \sum_{i=1}^{n-2}
\int_{t_{i-1}}^{t_{i}} \Biggl(\sum_{k=i+2}^{n} f_{k}^{t}(u) \int_{0}^{u}
f_{k}^{t}(s) \,dM_{s}\Biggr)\,dM_{u} \\
&=& n^{2H-1} \int_{0}^{t_{n-2}}
\Upsilon_{u,n}^{M} \,dM_{u},
\end{eqnarray*}
where
\[
\Upsilon_{u,n}^{M}=\sum_{k=i+2}^{n} f_{k}^{t}(u)\int_{0}^{u} f_{k}^{t}(s)
\,dM_{s},\qquad u \in[t_{i-1}, t_{i} ).
\]
%
%Put
%}(u)f^t_k(v)1_{(0,u)}(v).
%With this notation we have that
%S^{n} = n^{2H-1} \int_0^t \int_0^u \Upsilon_n(u,v) dM_vdM_u.
We use the following version of the Lenglart inequality: if $N$ is
a locally square integrable continuous martingale, then, for any
$\varepsilon>0$, $t>0$ and $A>0$,
%
%e3.20 ###
\begin{equation}\label{eq:1.18.11}
\PP\Bigl\{{\sup_{0 \leq s \leq t}}|N(s)| \geq\varepsilon\Bigr\} \leq
\frac{A}{\varepsilon^{2}}+\PP\{[N]_{t} \geq A\}.
\end{equation}
It follows from inequality \eqref{eq:1.18.11} that it is sufficient to
prove the
relation
%
%e3.21 ###
\begin{equation}\label{eq:1.18.12}
n^{4H-2}\int_{0}^{t_{n-2}} (\Upsilon_{u,n}^{M})^2\,d[M]{u} \stackrel
{\PP}
\to0,\qquad n \to\infty.
\end{equation}
First, using integration by parts, we estimate the function
\[
\Upsilon_{u,n}^{M}=\sum_{k=i+2}^{n} f_{k}^{t}(u) \biggl[f_{k}^{t}(u)M_{u}-
\int_{0}^{u} M_{s} (f_{k}^{t} (s) )'_s\,ds\biggr],\qquad u \in[t_{i-1},
t_{i} ).
\]
Clearly,
\[
(f_{k}^{t}(u) )_u' =
\biggl(\frac32-H \biggr)\int_{t_{k-1}}^{t_{k}} v^{H-1/2}
(v-u)^{H-5/2}\,dv.
\]
Therefore,
\begin{eqnarray}
|\Upsilon_{u,n}^{M}|&\leq& L \sum_{k=i+2}^{n}
(f_{k}^{t}(u))^{2}
\nonumber\\
&&{}+L \biggl(\frac32-H \biggr)\sum_{k=i+2}^{n}f_{k}^{t}(u)
\int_{0}^{u} \int_{t_{k-1}}^{t_{k}} v^{H-1/2}(v-s)^
{H-5/2}\,dv \,ds,\nonumber\\
\eqntext{u \in[t_{i-1}, t_{i} ).}
\end{eqnarray}

We estimate the terms separately: since $f_{k}^{t}(u)\leq
\frac{{t}^{H+1/2}}{n}
(t_{k-1}-u )^{H-3/2}$,
we have that, for $u \in[t_{i-1}, t_{i} )$,
\begin{eqnarray*}
\sum_{k=i+2}^{n}(f_{k}^{t}(u))^{2} &\leq& \frac{t^{2H+1}}{n^{2}}
\sum_{k=i+2}^{n} (t_{k-1}-u)^{2H-3} \\
&\leq& \frac{t^{2H+1}}{n^2}
(t_{i+1}-u)^{2H-3}
+ \frac{t^{2H+1}}{n} \int_{t_{i+1}}^{1}
(tx-u)^{2H-3} \,dx \\
&\leq& \frac{t^{4H-2}}{n^{2H-1}} + \frac{t^{2H}}{n}
\frac{(t_{i+1}-u)^{2H-2}}{2-2H}\\
&\leq& Cn^{1-2H}
\end{eqnarray*}
and
\begin{eqnarray*}
&&
\sum_{k=i+2}^{n}f_{k}^{t}(u) \int_{0}^{u} \int_{t_{k-1}}^{t_{k}}
v^{H-1/2}(v-s)^ {H-5/2}\,dv \,ds \\
&&\qquad\leq
C\sum_{k=i+2}^{n}f_{k}^{t}(u) \int_{t_{k-1}}^{t_{k}}
v^{H-1/2}(v-u)^ {H-3/2}\,dv
\\
&&\qquad\leq C \sum_{k=i+2}^{n}
(f_{k}^{t}(u))^{2}\leq Cn^{1-2H}.
\end{eqnarray*}

From these estimates, it follows that $
n^{4H-2}(\Upsilon_{u,n}^{M})^2\leq C $. Therefore, the
bounded majorant in \eqref{eq:1.18.12} exists. So, in order to establish
\eqref{eq:y-e-f-1}, it is sufficient to prove that
$\Upsilon_{u,n}^{M} n^{2H-1}\stackrel{\PP}\to0$, $0<u<t$. We have that
%
%e3.22 ###
\begin{eqnarray}\label{eq:dom1}\quad
&&
\E(\Upsilon_{u,n}^{M} n^{2H-1})^{2}\nonumber\\
&&\qquad= n^{4H-2}
\E\int_0^u \Biggl(\sum_{k=i+2}^{n} f_{k}^{t}(u)
f_{k}^{t}(s) \Biggr)^2 \,d [{M}]_{s},\\
\eqntext{u \in[t_{i-1},
t_{i} ).}
\end{eqnarray}

Similarly to previous estimates, we obtain that
\begin{eqnarray*}
&&n^{4H-2} \Biggl(\sum_{k=i+2}^{n}f_{k}^{t}(u)
f_{k}^{t}(s) \Biggr)^{2}\\
&&\qquad\leq C
n^{4H-2} \Biggl(\sum_{k=i+2}^{n}
\frac{1}{n^{2}} (t_{k-1}-u )^{H-3/2}
(t_{k-1}-s )^{H-3/2} \Biggr)^{2}\\
&&\qquad\leq C n^{4H-4} \biggl(\frac{1}{n}\sum_{k=i+2}^{n}
(t_{k-1}-u )^{2H-3} \biggr)^{2}
\\
&&\qquad\leq C n^{4H-4}
\biggl(\frac{n^{3-2H}}{n}+n^{2-2H} \biggr)^{2} \leq C \qquad\mbox{for
some } C>0.
\end{eqnarray*}
This means that the bounded
dominant in \eqref{eq:dom1} exists. Moreover,

%%%%%%%%%%%%%%%%%%%%%%%%%%%%%%%%%%%%%%%%%%%%%%%%%%%%%%%%%%%%%%%%%%%%%%%%%%%%%%%
%
\begin{eqnarray*}
&&
n^{2H-1}\sum_{k=i+2}^{n} f_{k}^{t}(u) f_{k}^{t}(s)\\
&&\qquad\leq
C n^{2H-1}\sum_{k=i+2}^{n} f_{k}^{t}(u) \cdot\frac{1}{n}
(u-s)^{H-3/2} \\
&&\qquad\leq C n^{2H-1} \cdot\frac{1}{n} \int
_{({i+1})/{n}}^{1} v^{H-1/2}(v-u)^{H-3/2}\,dv \cdot
(u-s)^{H-3/2}\rightarrow0
\end{eqnarray*}
for any $s<u$. Putting\vspace*{-3pt} together, this means that $\Upsilon_{u,n}^M
n^{2H-1}\stackrel{\PP}\to0$, $0<u<1$, whence $S^{n} \stackrel{\PP}
\rightarrow0$ and, consequently, \eqref{eq:y-e-f-1} holds.
Next, consider the sums
\[
n^{2H-1} \sum_{k=n {s/t} +2}^n \int_{t_{k-2}}^{t_{k-1}}f^t_k(u)
\int
_{t_{k-2}}^u f^t_k(v)\,dM_v \,dM_u
\]
and
\[
n^{2H-1} \sum_{k=n {s/t} +1}^n \int_{t_{k-1}}^{t_k}g^t_k(u)
\int
_{t_{k-1}}^u g^t_k(v)\,dM_v \,dM_u.
\]
The assumptions of Lemma \ref%
{le:img} are satisfied with martingales
\[
N^{n,k}_v := n^{H-1/2} \int_{t_{k-2}\wedge v}^{t_{k-1}\wedge v} f^t
_k(u)\,dM_u
\]
and
\[
\tilde N^{n,k}_v := n^{H-1/2} \int_{t_{k-1}\wedge v}^{t_k\wedge
v} g^t
_k(u)\,dM_u.
\]
Indeed, property (ii) follows from \eqref{eq:f-y-12} and \eqref{eq:y3},
and property (iii)
can be easily checked.
Hence, both sums are of the order $o_\PP(1)$.
The next statement is an immediate consequence of Lemma \ref{le:y1} and
\eqref{eq:y-e-f-1}.
There exist two constants $C_{1}>0$, $C_{2}>0$ such that
%
%
%e3.23 ###
\begin{eqnarray}\label{eq:gen}
C_{1}t^{2H-1}\int_s^t u^{2H-1}\,d[M]_u
&\leq&{\PP}-\lim_{{n}\rightarrow\infty}
{n}^{2H-1}
\sum_{k={{n}s}/{t}+1}^{{n}}
(I^{n,1}_k)^{2}\nonumber\\[-8pt]\\[-8pt]
&\leq& C_{2}t^{4H-2}([M]_{t}-
[M]_{s}).\nonumber
\end{eqnarray}

Similarly, one can show that the cross product sums with $%
i\neq j$ satisfy $n^{2H-1} %\times\break
\sum_k I^{n,i}_kI^{n,j}_k = o_\PP(1) $. Indeed,
let $i=1$ and $j=2$; other cases can be considered similarly. We have
that, in this case,
\[
{n}^{4H-2} \E\Biggl(\sum_{k=1}^{n}I_{k}^{n,1}I_{k}^{n,2} \Biggr)^{2}
%{n}^{4\alpha} E
%(\sum_{k=1}^{n}(I^{k}_{1})^{2}(I_{k}^{2}))^{2}
=n^{4H-2}\E\sum_{k=1}^n(I_{k}^{n,1})^2J_k^{n,2},
\]
where
$J^{n,2}_k=\int^{t_{k-1}}_{t_{k-2}}(f_k^t(s))^2\,d[M]_s$, since
$I_{k}^{n,1}$, $I_{k}^{n,2}$, $I_{k}^{n,3}$
are pairwise orthogonal. Moreover, the product sum $n^{2H-1} \sum_k
I^{n,i}_kI^{n,j}_k$ can be considered as a final value of a square
integrable martingale with quadratic characteristic $\sum
_{k=1}^n(I_{k}^{n,1})^2J_k^{n,2}$. So, it follows from the Lenglart
inequality that it is sufficient to prove the relation
%
%e3.24 ###
\begin{equation}\label{eqlim}{n}^{4H-2} \sum
_{k=1}^{n}(I_{k}^{n,1})^{2}J_k^{n,2}\stackrel{\PP}\to0.
\end{equation}
According to \eqref{eq:gen}, we have that
\[
\PP-\lim_{n\to\infty}n^{2H-1}\sum_{k=1}^{n}(I_{k}^{n,1})^2
\leq C_2t^{4H-2} [M]_{t}
\]
and, also,
\[
n^{2H-1}\max_{1\leq k\leq n}\int
^{t_{k-1}}_{t_{k-2}}(f_k^t(s))^2\,d[{M}]_s\leq
\biggl(H-\frac12 \biggr)^{-2}\max_{1\leq k\leq n}([{M}]_{t_
{k-1}}-[{M}]_{t_{k-2}})\stackrel{\PP} \rightarrow0,
\]
whence
\eqref{eqlim} follows.
%$$(n^{2H-3}t^{2H+1}\sum_{i=k+2}^{n} (f_i^t(u))^{2}+(g_k^t(u))^{2})1_{
% + n^{2H-1}(f_k^t(u))^{2}1_{\{u\in[t_{k-2}, t_{k-1})\}}
%$$

We are now ready to
finish the proof of Theorem \ref{t:levy} in the case $H>\frac12$.

%s3.4 ###
\subsection{Completion of the proof for the case $H >
\frac12$}
\label{ss:end-of-proof}

Suppose, for the moment, that we consider the
fixed interval $[0,t]$. By using our estimates, we can conclude that
for rational $s$,
consequently for any $s<t$, the following claims hold:

(d) there exist two constants, $C_{1}>0$ and $C_{2}>0$, such that
\[
C_{1}\int_{s}^{t}u^{2H-1}\,d[M]_{u}\leq
t-s \leq C_{2}t^{2H-1} ([M]_{t}-
[M]_{s});
\]
this estimate can be rewritten in
terms of $W$ and $[W]$ (recall that $W_t = \int_0^t s^{H-1/2}\,dM_s
$) as
\[
C_{1}([W]_{t}-[W]_{s})\leq t-s
\leq C_{2}t^{2H-1} \int_{s}^{t}u^{1-2H}\,d[W
]_{u};
\]

(e)
\[
\PP-\lim_{n\rightarrow\infty} n^{2H-1}
\sum_{k=n{s/t}+1}^{n} (\triangle
X_{t_{k}})^{2} =\PP-\lim_{n\rightarrow\infty} \int_{s}^{t}
\varphi^{t}_{n}(u)\,d[M]_{u},
\]
where we can take $\varphi^{t}_{n}(u)$ from \eqref{eq:est-2f}, \eqref
{eq:f-y-12} and \eqref{eq:y3}, and they equal
\begin{eqnarray*}
\varphi^{t}_{n}(u)&=&\Biggl(n^{2H-3}t^{2H+1} \sum_{i=k+2}^{n}
(f_i^t(u))^{2}+n^{2H-1}(g_k^t(u))^{2}\Biggr)1_{\{u\in[t_{k-1},t_k)\}}
\\
&&{}+n^{2H-1}(f_k^t(u))^{2}1_{\{u\in[t_{k-2},t_{k-1})\}}.
\end{eqnarray*}
Clearly, $\varphi^{t}_{n}(u)$ are positive, bounded, nonrandom
functions and it follows from \eqref{eq:f-down} that they are separated
from $0$ by some
constant multiplied by $u^{2H-1}$.

From the left-hand side of (d), it follows that $[W]_{t}$ is absolutely
continuous with respect to the Lebesgue measure,
so $[W]_{t}=\int_{0}^{t}\theta_{s}\,ds$, where
$\theta_{s}$ is a bounded, possibly random, variable. From the
right-hand side of (d), it follows that
\[
\int_{s}^{t} u^{1-2H}\theta_{u}\,du\geq\frac{1}{C_{2}}
(t^{2-2H}-st^{1-2H})\geq C_{3}
(t^{2-2H}-s^{2-2H})= C_{3}\int_{s}^{t} u^{1-2H}\,du.
\]
This means that
\[
\int_{s}^{t} u^{1-2H} (\theta_{u}-C_{3})\,du\geq0,
\]
%
%Evidently, for any set $A\in\mathcal{F}$
%$$\int_{A}\int_{s}^{t} u^{2H-1} (\theta_{u}-C_{3})\,du dP\geq0.$$
%Now, let the set $D \in\sigma\{ \F\times\mathcal B [\delta,T]\}$,
%and let $\delta>0$ be fixed. Then $\mu(D)<\infty$, where $\mu= P
%
%%%%%%%%%%%%%%%%%%%%%%%%%%%%%%%%%%%%%%%%%%%%%%%%%%%%%%%%%%%%%%%%%%%%%%%%%%%%%%%%
%
%By the theorem of approximation of measurable sets, for any
%$\varepsilon>0$ there exists a collection of the sets
%
%$$\{D_{i}=B_{i}\times[s_{i},t_{i}],
% B_{i}\in\F, [s_{i},t_{i}] \in\mathcal B[\delta,T]\},$$
%such that
%$$\mu( (D \setminus\bigcup_{i=1}^{k} D_{i})\bigcup
% (\bigcup_{i=1}^{k}D_{i} \setminus D ) )<\varepsilon.$$
%Therefore, since $u^{2H-1}(\theta_{u}-C_{3})$ is bounded on $D$, we
%get that
%$\int_{D} u^{2H-1}(\theta_{u}-C_{3})d\mu\geq0.$
%
%Now, set $D= \{ (\omega, u): \theta_{u}-C_{3}<0,\text{ and }
%u\geq\delta\}$ and we immediately obtain that $\mu(D)=0$.
whence we immediately obtain that $\theta_u(\omega)>C_3>0$ for almost
all $u, \omega$, concluding that $[W]$ is equivalent to the Lebesgue
measure and so $ W_{t}=\int_{0}^{t}\theta_{s}^{{1/2}}\,dV_{s}$,
where $\{V_{s}, { \F}_{s}, s \geq0 \}$ is some Wiener process.

Now, if we perform all of the same calculations as before, but for ``true''
fractional Brownian motion $B^H_{t}$, we obtain that
\begin{eqnarray*}
\PP-\lim_{n\rightarrow\infty}
n^{2H-1}\sum_{k=n{s/t}+1}^{n}(\triangle
B^H_{t_{k,n}})^{2}&=&\PP-\lim_{n\rightarrow\infty}\int
_{s}^{t}\varphi
_{s}^{n}s^{2H-1}\,ds
\\ &=&t ^{2H-1}(t-s).
\end{eqnarray*}
(It is sufficient to take $s=0$.) Therefore,
$\PP-\lim_{n\rightarrow\infty}\int_{s}^{t}\psi_{u}^{n}\,du=0$, where
$\psi_{u}^{n}=u^{2H-1}\varphi_{u}^{n}(\theta_{u}-1)$.

%Consider any set $D\in\sigma\{ \F\times\mathcal B [\delta,T]\}$
From this, we obtain that $\theta_{u}\equiv1$ [otherwise, consider
the set $D=\{(\omega, u)\dvtx\theta_{u}>1+\alpha, {\mbox{ or }}
\theta_{u}<1-\alpha\}$ for $\alpha>0$; clearly, it has zero
measure].

%s4 ###
\section{\texorpdfstring{Proof of Theorem \protect\ref{t:levy}:
$H<\frac12$}{Proof of Theorem 1.1: $H<\frac12$}}

For $H<\frac12$, we use, in general, principally the same ideas.
However, technical details are different. Indeed, it is well known
(see, e.g., \cite{nvv}) that the kernel $z(t,s)$ participating in the
representation of $X$ via $M$ or $W$ [see \eqref{eq:xw}] is more
complicated in the case $H<\frac12$. The brackets of the martingales
that are to be estimated as before also have an additional singularity
because the power $2H-1$, or any other power of such a form, is now
negative. Therefore, the proofs are more technical and the reasons for
this will be mentioned below in all relevant places.

%s4.1 ###
\subsection{Starting point}

At first, consider the H\"{o}lder properties of the processes
involved. We can note the following: since $H<\frac{1}{2}$,
it is very simple to prove, using integration by parts, that the
process $Y$ has the same
H\"{o}lder properties as $X$, that is, it is H\"{o}lder up to order
$H$. Further, it follows from Lemma 2.1 \cite{nvv}
that $M$ is H\"{o}lder up to order $\frac{1}{2}$. Therefore, for any
$0<s_0\leq s<t\leq T$ and $\beta<\frac{1}{2}$, there exists a constant
$K=K_{s_0, \beta}$ such that $|W_t-W_s|\leq K_{s_0,\beta}(t-s)^\beta$.
Now, it is more convenient to consider $W$ instead of $M$. We shall
show the inequality
%
%e4.1 ###
\begin{equation} \label{eq:tochtonado}
C_1 ( [W]_t - [W]_s ) \le t-s \le C_2 ( [W]_t - [W]_s
)
\end{equation}
first for arbitrary $t> 0$ and $s \in\mathcal{R}_t$, $s < t $. Recall
that we can assume the processes $%
W$ and $[W] $ to be bounded, as in Section \ref{ss:basic}.

For $H<\frac12$, we use the following representation result,
which can be proven as \cite{nvv}, Theorem 5.2.
\begin{lemma}
\label{le:xw} Assume that $H< \frac12 $ and that
properties \textup{(a)} and
\textup{(c)} hold. The process $X$ then has the representation
%
%e4.2 ###
\begin{equation} \label{eq:xw}
X_t =\int_{0}^{t}{z(t, s) \,dW_s}
\end{equation}
with the kernel
\begin{eqnarray*}
z(t, s) &=& \biggl(\frac{s}{t} \biggr)^{1/2-H}(t-s)^{H-1/2} \\
&&{} - (H-1/2)s^{1/2-H}
\int_s^t{u^{H-3/2}(u-s)^{H-1/2}\,du}.
\end{eqnarray*}
\end{lemma}

Put
\[
p^t_k (z) = \int_{t_{k-1}}^{t_k} \biggl(\frac{z}{u}\biggr)^{1/2 -H} (u-z)
^{H-3/2} \,du
\]
for $z< t_{k-1}$.

Using Lemma \ref{le:xw} and integration by parts, we can now
write the increment of $X$ as
\begin{eqnarray*}
X_{t_k} - X_{t_{k-1}} & = & \biggl(\frac{1}{2}-H \biggr) \int_0^{t_{k-2}}
p^t_k(s)\,dW_s\\
& &{} + \biggl(\frac{1}{2}-H \biggr) \int_{t_{k-2}}^{t _{k-1}} p^t_k(s)
\,dW_s \\
& &{} + \int_{t_{k-1}}^{t_k} \biggl( \frac{s}{t_k} \biggr)^{1/2-H} (
t_k-s )^{H-1/2}\,dW_s \\
& &{} + \biggl(\frac{1}{2}-H \biggr) \int_{t_{k-1}}^{t_k} s^{1/2-H}
\int_s^{t_k}u^{H-3/2} (u-s)^{H-1/2}\,du \,dW_s \\
& =&\!: J_k^{n,1} + J_k^{n,2}+ J_k^{n,3} + J_k^{n,4}.
\end{eqnarray*}

Clearly,
\begin{eqnarray*}
&&\lim_{n \to\infty}n^{2H-1}
\sum_{k=n{s/t}+2}^n ( \Delta X_{t_{k}} )^2 \\
&&\qquad = \lim_{n
\to\infty} n^{2H-1} \Biggl( \sum_{k = n{s/t}+2}^n (
J_k^{n,1} )^2
+ \sum_{k = n{s/t}+2}^n (
J_k^{n,2} + J_k^{n,3} + J_k^{n,4} )^2\\
& &\qquad\quad\hspace*{106.7pt}{} +
2\sum_{k = n{s/t}+2}^n J_k^{n,1} (J_k^{n,2}+ J_k^{n,3} +
J_k^{n,4} ) \Biggr).
\end{eqnarray*}
As before,
%
%e4.3 ###
\begin{equation} \label{eq:vottak}\lim_{n \to\infty}n^{2H-1}
\sum_{k=n{s/t}+2}^n ( \Delta X_{t_{k,n}} )^2 \stackrel
{L^1(\PP)}\to
t^{2H-1}(t-s).
\end{equation}
First, estimate
\[
\lim_{n \to\infty} n^{2H-1} \sum_{k
= n{s/t}+2}^n (J_k^{n,1} )^2
\]
from below and
above.
We start with the analog of Lemma \ref{le:y1}.

%s4.2 ###
\subsection{Two-sided estimates for the sums $n^{2H-1}
\sum_{k ={ns/t}+2}^n \int_0^{t_{k-2}} ( p_k^t (z) )^2\, d
[ W]_z$ and $n^{2H-1} \sum_{k= n{s/t}+2}^n (J_k^{n,1} )^2$}
\label{ss:upper-p-e}

Put
\[
j^{n,1}= n^{2H-1} \sum_{k={ns/t}+2}^n \int_0^{t_{k-2}} ( p_k^t
(z) )^2 \,d [ W ]_z.
\]
We decompose this sum as in the case of the proof for $H>\frac12 $ [see
\eqref{eq:est-2} and \eqref{eq:est-2a}]:
\begin{eqnarray*}
j^{n,1}:= n^{2H-1} \Biggl( \sum_{i=1}^{n {s/t}} \sum_{k =n
{s/t}%
+2}^n + \sum_{i =n {s/t}+1}^{n-2} \sum_{k=i+2}^n \Biggr)
\int_{t_{i-1}}^{t_i} ( p_k^t (u) )^2 \,d [ W]_u.
\end{eqnarray*}
Clearly, for $s\leq t_{k-2}$,
\[
p_k^t (s) \leq\biggl( ( t_{k-1}-s )^{H-3/2}
\frac{t}{n} \biggr)\wedge\biggl(\frac{1}{1/2-H} \biggl(\frac
tn \biggr)^{H-1/2} \biggr).
\]
Therefore, for $n$ such that $n\frac st\in\N$, we have that
%
%
%e4.4 ###
\begin{eqnarray}\label{111}
j^{n,1} &\leq& n^{2H-2} t\int_0^s \biggl( s + \frac{t}{n}
-u \biggr)^{2H-2} \,d[W]_u \nonumber\\
&&{}+n^{2H-3} t^2\int_0^s \biggl( s +
\frac{t}{n} -u \biggr)^{2H-3} \,d[ W]_u \nonumber\\[-8pt]\\[-8pt]
&&{}+\frac{t}{2-2H} \biggl(
\frac{t}{n} \biggr)^{2H-2} n ^{2H-2}\sum_{i = n{s/t}+1}^{n-2}
\int_{t_{i-1}}^{t_{i}} \,d[W]_u \nonumber\\
&&{}+t^2 n^{2H-3}\sum_{i=n{s/t}+1}^{n-2} \int_{t_{i-1}}^{t_{i}}\, d[W]_u
\biggl(\frac{t}{n} \biggr)^{2H-3}.\nonumber
\end{eqnarray}

We divide the integral $\int_0^s (s + \frac{t}{n} -u )^{2H-2}\,
d[W]_u$ into two parts, $\int_0^{s/2} (s + \frac{t}{n} -u
)^{2H-2}\,
d[W]_u$ and $\int_{s/2}^s (s + \frac{t}{n} -u )^{2H-2}\,
d[W]_u$. The first integral can be estimated as
\[
\int_0^{s/2} \biggl(s + \frac{t}{n} -u \biggr)^{2H-2}
\,d[W]_u\leq\biggl(\frac{s}{2} + \frac{t}{n} \biggr)^{2H-2}[W]_{{s/2}},
\]
whence $n^{2H-2} t\int_0^{s/2} ( s + \frac{t}{n}
-u )^{2H-2} \,d[W]_u\rightarrow0$ as $n\rightarrow\infty$ a.s. As
for the second part, we apply the following inequality from \cite{nvv}:
let the function $f\dvtx[a,b]\rightarrow R$ be H\"{o}lder on $[a,b]$ of
order $\beta$, $|f(t)-f(s)|\leq K|t-s|^\beta$. Then, for any $\rho
>-1+\beta$ and $b<v$, we have that
%
%e4.5 ###
\begin{equation}\label{eq:4.24}
\quad
\biggl|\int_a^b (v-u)^\rho\, df(u)\biggr|\leq K
\biggl(1+ \biggl|\frac{\rho}{\rho+\beta} \biggr| \biggr)\bigl((v-b)^{\rho+\beta}+(v-a)^{\rho
+\beta}\bigr).
\end{equation}
According to the H\"{o}lder properties of $W$ mentioned above, we can
take any $0<\beta<\frac{1}{2}$ and define, for any $r\in(\frac
{s}{2},t]$, the random variable
\[
K_r(\omega)=\sup_{{s/2}\leq u<v\leq r}\frac
{|W_v-W_u|}{(v-u)^\beta}.\vadjust{\goodbreak}
\]
Clearly, $\PP\{K_t(\omega)\geq N\} \to0$ as $N\rightarrow\infty$.
Therefore, it is enough to prove that
$\int_{{s/2}}^{s\wedge\tau_N} (s + \frac{t}{n} -u
)^{2H-2}\,d[W]_un^{2H-2}\stackrel{\PP} \to0$ as $N\rightarrow\infty$ for
any $N>1$, where $\tau_N=\inf\{r\geq\frac{s}{2}\dvtx K_r\geq N\}\wedge t$.
According to the Burkholder--Gundy inequality and \eqref{eq:4.24},
\begin{eqnarray}
&&
n^{2H-2}\E\biggl(\int_{{s/2}}^{s\wedge\tau_N} \biggl(s + \frac
{t}{n} -u \biggr)^{2H-2}\,d[W]_u \biggr)
\nonumber\\
&&\qquad
\leq Cn^{2H-2}\E\biggl(\int_{{s/2}}^{s\wedge\tau_N} \biggl(s + \frac
{t}{n} -u \biggr)^{H-1}\,dW_u \biggr)^2
\nonumber\\
&&\qquad
\leq CN^2n^{2H-2} \biggl(\frac{\beta}{H+\beta-1} \biggr)\nonumber\\
&&\qquad\quad{}\times \biggl( \biggl(\frac
{t}{n} \biggr)^{H+\beta-1}+ \biggl(\frac{s}{2}+\frac{t}{n} \biggr)^{H+\beta
-1} \biggr)^2\rightarrow0 \nonumber\\
\eqntext{\mbox{as } n\rightarrow\infty.}
\end{eqnarray}
Finally, we obtain that $n^{2H-2} \int_0^s ( s + \frac{t}{n}
-u )^{2H-2} \,d[W]_u\stackrel{\PP} \to0$ as $n\rightarrow\infty$.
%
%Lemma 2.1 \cite{nvv}, can be
%estimated as $$ |\int_0^s ( s + \frac{t}{n}
%-u )^{2H-2} d[W]_u | \leq C(\omega)
% (s + \frac{t}{n} -s )^{2H-2+\beta},
%$$ for some random variable $0<C(\omega)<\infty$, where $\beta$
%is H\"{o}lder index of $[W]_u$. Evidently, $\beta> 0$,
%and it holds that
%$$ \int_0^s ( s + \frac{t}{n} -u )^{2H-2} d[
%W]_u\cdot n^{2H-2}\sim n^{2H-2} (

The same is true for
\[
\int_0^s \biggl( s + \frac{t}{n} -u \biggr)^{2H-3} \,d[W]_u\cdot n^{2H-3}.
\]
The last two integrals from \eqref{111} admit
the obvious estimate $t^{2H-1} C_2 ([W]_t - [W]_s )$.

The ``remainder'' term for $\sum(
J^k_1 )^2$, that is, the difference between $\sum(
J^k_1 )^2$ and $j^{n,1}$, equals
\begin{eqnarray*}
R_n &:=& n^{2H-1}
\sum_{k=n{s/t}+2}^{n} \int_0^{t_{k-2}} \biggl( \int_0^z
p_k^t (v) \,dW_v \biggr)\\
&&\hspace*{92.8pt}{}\times p_k^t (u)\,dW_u.
\end{eqnarray*}

For technical simplicity, it is enough to consider $\sum_{k=3}^{nr}$
for any $r \in\N$, instead of
$\sum_{k=n{s/t}+2}^{n} = -\sum^{n{s/t}+1}_{k=3}+
\sum_{k=3}^n$. We obtain that
\begin{eqnarray*}
\E(R_n)^2 &=& n^{4H-2}
\E\Biggl( \sum_{k=3}^{nr} \sum_{i=1}^{k-2}
\int_{t_{i-1}}^{t_{i}} \int_0^u p_k^t (v) \,d {W}_v
\cdot p_k^t (u) \,d {W}_u \Biggr)^2
\\
&=& n^{4H-2} \E\Biggl( \sum_{i=1}^{nr-2}
\sum_{k=i+3}^{nr} \int_{t_{i-1}}^{t_{i}} \int_0^u p_k^t
(v) \,d{W}_v \cdot p_k^t (u) \,d{W}_u \Biggr)^2
\\
&=& n^{4H-2}\sum_{i=1}^{nr-2} \E\int_{t
_{i-1}}^{t_{i}} \Biggl(\sum_{k=i+3}^{nr} \int_0^u p_k^t
(v) \,d{W}_v \cdot p_k^t (u) \Biggr)^2 \,d[W]_u.
\end{eqnarray*}
Let us estimate
\begin{eqnarray*}
\biggl| \int_0^u p_k^t (v) \,dW_v \biggr|&=&
\biggl|p_k^t (u) W_u - \int_0^u {W}_v (
p_k^t (v) )'_v\,dv \biggr|\\
&\leq& L | p_k^t
(u) | + L \biggl| \int_0^u (p_k^t (v) )'_v \,dv
\biggr|.
\end{eqnarray*}
We have that %$z \leq t_{k-2}$
%$$
%p_k^t (u) \leq\int_{t_{k-1}}^{t_{k}}
%(u-z)^{\alpha-1} \,du \leq(t_{k-1}-z )^{\alpha}
%$$
%
\[
\biggl|\int_0^u (p_k^t (v) )'_v \,dv \biggr| = |
p_k^t (u)-p_k^t (0) | \leq C \biggl(
\frac{t}{n} \biggr)^{H-1/2} \qquad\mbox{for some } C>0.
\]

Moreover,
\begin{eqnarray*}
n^{2H-1} \Biggl( \sum_{k=i+3}^{nr} p_k^t
(u) \Biggr)^2 &\leq& n^{2H-1} \biggl( \int_{t_{i+1}}^{tr}
(v-u)^{H-3/2} \,dv \biggr)^2
\\
&=& Cn^{2H-1} [ -
(tr-u)^{H-1/2} + ( t_{i+1}-u )^{H-1/2} ]^2 \\
&\leq& C
\end{eqnarray*}
and the integrand
\[
n^{4H-2} \Biggl( \sum_{k=i+2}^{nr}
\int_0^u p_k^t (v) \,d {W}_v \cdot p_k^t(u)
\Biggr)^2 \leq C,
\]
that is, the integrable dominant exists.
Therefore, it is sufficient to establish that for any $u$,
\[
n^{2H-1} \sum_{k=i+3}^{nr} \int_0^u p_k^t (v) \,
d{W}_v \cdot p_k^t(u) \stackrel{\PP} \to0.
\]
We take the
mathematical expectation in the left-hand side and obtain that
\[
n^{4H-2} \E\int_0^u \Biggl(
\sum_{k=i+3}^{nr} p_k^t (v) p_k^t (u) \Biggr)^2 \,d[W]_v.
\]
Also, here, the bounded dominant exists. Indeed,
\[
n^{4H-2} \Biggl( \sum_{k=i+3}^{nr} p_k^t (v) p_k^t
(u) \Biggr)^2 \leq n^{2H-1} \Biggl( \sum_{k=i+2}^{nr} p_k^t
(v) \Biggr)^2 \leq C,
\]
as before. Further, we must prove that
\[
n^{2H-1} \sum_{k=i+3}^{nr} p_k^t (v) p_k^t (u) \to0
\]
for all fixed $0 < v < u$. We have that
\begin{eqnarray}
&&
n^{2H-1}
\sum_{k=i+2}^{nr} p_k^t (v) p_k^t (u) \nonumber\\
&&\qquad\leq n^{2H-1}
\sum_{k=i+3}^{nr} \int_{t_{k-1}}^{t_{k}} (s-u)^{H-3/2} \,ds
\int_{t_{k-1}}^{t_{k}} (s-v)^{H-3/2} \,ds\nonumber\\
&&\qquad\leq n^{2H-1} \sum_{k=i+3}^{nr} (
t_{k-1}-u )^{H-3/2} \frac{1}{n} \int_{t_{k-1}}^{t_{k}}
(s-v)^{H-3/2} \,ds
\nonumber\\
&&\qquad\leq n^{2H-2} ( t_{i+2}-u )^{H-3/2}
\int_{t_{i+2}}^{tr} (s-v)^{H-3/2} \,ds
\nonumber\\
&&\qquad\leq C n^{H-3/2} (u-v)^{H-3/2} \to0 \nonumber\\
\eqntext{\mbox{as } n \to
\infty\mbox{ for any } 0<v<u.}
\end{eqnarray}
From all of these estimates, the remainder term $R_n \stackrel{\PP}
\to
0$.

For the lower bounds, we return to $[
M]$ instead of $[W]$:
\begin{eqnarray*}
j^{n,1} &=& n^{2H-1}
\sum_{k = n{s/t}+2}^{n} \int_0^{t_{k-2}} (f_t^k (u)
)^2 \,d[M]_u
\\
&\geq& t^2 n^{2H-3}\sum_{k =
n{s/t}}^{n} ( t_{k} )^{2H-1} \int_0^{t_{k-2}}
(t_{k}-u )^{2H-3} \,d[M]_u
\\
&\geq& t^2 n^{2H-3} \Biggl( \sum_{i=1}^{n {s/t}-1}
\sum_{k=n{s/t}+2}^n + \sum_{i=n{s/t}+1}^{n-2}
\sum_{k=i+2}^n \Biggr) (t_k)^{2H-1} \\
&&{}\times\int_{t_{i-1}}^{t_{i}} (
t_{k}-u )^{2H-3} \,d[M]_u
\\
&=& C t^{2H+1} n^{2H-2}
\sum_{i=n{s/t}+1}^{n-2} \int_{t_{i-1}}^{t_{i}} \frac{1}{t}
\bigl( ( t_{i+2}-u )^{2H-2}- (t-u)^{2H-2} \bigr)\,
d[M]_u.
\end{eqnarray*}
Note that
\begin{eqnarray*}
&&
n^{2H-2}
\sum_{i=n{s/t}+1}^{n-2} \int_{t_{i-1}}^{t_{i}}
(t-u)^{2H-2}\,d[M]_u
\\
%
%$$n^{2\alpha-1} \int_{s}^{t-2/n} (t-z)^{2\alpha-1} d [M\rangle_z
%
&&\qquad
\sim\biggl(t-t+\frac{2}{n} \biggr)^{2H-2+\beta} \cdot
n^{2H-2} \to0\qquad \mbox{as } n \to\infty.
\end{eqnarray*}
Therefore,
\begin{eqnarray*}
&&\lim_{n
\to\infty} n^{2H-1} \sum_{k=n{s/t}+1}^{n} (
J^k_1 )^2
\\[-1.5pt]
&&\qquad\geq C t^{2H} n^{2H-2}
\sum_{i=n{s/t}+1}^{n-2} \int_{t_{i-1}}^{t_{i}} (
t_{i+2}-u )^{2H-2}\,d[M]_u
\\[-1.5pt]
&&\qquad\geq C t^{2H} n^{2H-2} \sum_{i=n{s/t}+1}^{n-2}
( t_{i+2}-t_{i-1} )^{2H-2} \int_{t _{i-1}}^{t_{i}}
d[M]_u.
\end{eqnarray*}

Combining this with the upper estimate
and taking into account the estimate of the remainder term, we have
%
%e4.6 ###
\begin{eqnarray} \label{eq:half-a-3}
C_1 t^{4H-2}
([M]_t - [M]_s )
&\leq& \lim_{n
\to\infty} n^{2H-1} \sum_{k=n{s/t}+2}^{n} (
J^k_1 )^2\nonumber\\[-8.5pt]\\[-8.5pt]
&\leq& C_2 t^{2H-1} ([W]_t -
[W]_s ).\nonumber
\end{eqnarray}
[Note that, for $H\in(1/2,1)$, we
have obtained opposite estimates.] Also, note that we cannot
immediately estimate
$\sum( J^k_i )^2$, $i>1$, from above. Indeed, the
integrand of the form $ ( t \frac{l}{n}-u )^{H-1/2}$ that
admits the estimate $< (\frac{1}{n} )^{H-1/2}\to0$ for
$H\in(1/2,1)$, now, for $H\in(0, 1/2)$, tends to $\infty$. So, we can
mention that $\sum_{k=n{s/t}+2}^{n} ( J^k_2 + J^k_3 +
J^k_4 )^2 \geq0$, intend to prove that $\sum J^k_1 ( J^k_2
+ J^k_3
+ J^k_4 ) \to0$, and, from this, condition
(b) [or \eqref{eq:vottak}] and \eqref{eq:half-a-3}, obtain the
following estimate from above:
\[
C_1
t^{2H-1} ([M]_t - [M]_s )
\leq(t-s).
\]
In the sequel, we realize this plan.

%s4.3 ###
\subsection{Auxiliary estimates for ``mixed'' terms}

We will show that as $n\to\infty$, we have
%
%e4.7 ###
\begin{equation} \label{eq:ac-2}
n^{2H-1}\sum_k J^{n,1}_k (J^{n,2}_k + J^{n,3}_k + J^{n,4}_k )
\stackrel{\PP}{\to} 0.
\end{equation}
It is sufficient to estimate the sums from $k=2$ up to $k=n$. By
applying the Lenglart inequality to $n^{2H-1} \sum_{k=2}^{n}
J^{n,1}_kJ^{n,2}_k $ as well as to the final value of corresponding
martingale, we obtain that it is sufficient to prove that
\begin{eqnarray*}
&&n^{4H-2} \sum_{k=2}^{n} \biggl( \int_0^{t_{k-2}}
\int_{t_{k-1}}^{t_{k}} \biggl(
\frac{s}{u} \biggr)^{1/2-H}(u-s)^{H-3/2} \,du\,
d{W}_s \biggr)^2\\[-1.5pt]
&&\quad{}
\times\biggl( \int_{t_{k-2}}^{t_{k-1}}
\biggl(\int_{t_{k-1}}^{t_{k}} \biggl(
\frac{s}{u} \biggr)^{1/2-H}(u-s)^{H-3/2} \,du \biggr)^2 \,
d[W]_s \biggr)\\[-1.5pt]
&&\qquad
\leq C n^{4H-2} \sum_{k=2}^{n}
\biggl( \int_0^{t_{k-2}} p_k^t (s) \,d {W}_s \biggr)^2
\int_{t_{k-2}}^{t_{k-1}} ( t_{k-1}-s )^{2H-1} \,d[
W]_s \stackrel{\PP} \to0.
\end{eqnarray*}

Integrate the last integral by parts:
\begin{eqnarray*}
&&
\int_{t_{k-2}}^{t_{k-1}}
( t_{k-1}-s )^{2H-1} d[W]_s \\
&&\qquad=
( t_{k-1}- t_{k-2} )^{2H-1} ([
W]_{t_{k-1}} - [W]_{t_{k-2}} )
\\[-1pt]
&&\qquad\quad{}
- (2H-1)\int_{t_{k-2}}^{t_{k-1}} ( t_{k-1}-s )^{2H-2} (
[W]_{t_{k-1}} - [W]_{s} ) \,ds
\\[-1pt]
&&\qquad\leq C n^{1-2H} \Delta
[W]_{t_{k-1}}
 + C \int_{t_{k-2}}^{t_{k-1}}
( t_{k-1}-s )^{2H-2} ( [W]_{t_{k-1}} - [
W]_{s} ) \,ds.
\end{eqnarray*}

Now, recall that
\begin{eqnarray*}
\biggl( \int_0^{t_{k-2}} p_k^t (s) \,dW_s \biggr)^2
&=& \int_0^{t_{k-2}} ( p_k^t
(s) )^2 \,d[W]_s \\
&&{}+ 2\int_0^{t_{k-2}}
\int_0^s p_k^t (v) \,d W_v\cdot p_k^t (s)
\,d{W}_s.
\end{eqnarray*}
Clearly,
\[
\sigma^{n,1} :=
n^{2H-1} \sum_{k=2}^{n} \int_0^{t_{k-2}} ( p_k^t
(s) )^2 \,d[W]_s \leq j^{n,1},
\]
so, it is bounded in
probability and, similarly to $R_n$,
\[
\sigma^{n,2} := n^{2H-1} \sum_{k=2}^{n}
\int_0^{t_{k-2}} \int_0^s p_k^t (v) \,d {W}_v \cdot p_k^t
(s) \,d {W}_s \stackrel{\PP} \to0 \qquad\mbox{as } n \to\infty.
\]

Therefore,
\begin{eqnarray*}
&&
n^{4H-2} \sum_{k=2}^{n} \biggl( \int_0^{t_{k-2}}
p_k^t (s) \,d {W}_s \biggr)^2 \cdot C n^{1-2H}\Delta
[W]_{t_{k-1}}
\\
&&\qquad\leq C\sigma^{n,1} \cdot\max_k \Delta[
W]_{t_{k-1}} + C\sigma^{n,2} \cdot\max_k \Delta
[W]_{t_{k-1}} \stackrel{\PP} \to0,\qquad n \to
\infty.
\end{eqnarray*}

Also,
\begin{eqnarray*}
&&
n^{4H-2} \sum_{k=2}^{n} \biggl( \int_0^{t_{k-2}}
p_k^t (s) \,d {W}_s \biggr)^2 \cdot
\int_{t_{k-2}}^{t_{k-1}} ( t_{k-1}-s )^{2H-2}
\\
%
%$$
%n^{2\alpha}
% (\sigma_1^n +\sigma_2^n ) \int_{t_{k-2}}^{t
%(\omega) (t_{k-1} -s )^{1-\varepsilon}\,ds \leq
%$$
%
&&\qquad\quad\hspace*{132.9pt}{}
\times([{W}]_{t_{k-1}} -
[W]_{s} ) \,ds \\
&&\qquad\leq C
(\omega) (\sigma^{n,1} +\sigma^{n,2} ) n^{2H-1}
\int_{t_{k-2}}^{t_{k-1}} (
t_{k-1}-s )^{2H-1-\varepsilon} \,ds
\\
&&\qquad\leq
C(\omega) (\sigma^{n,1} +\sigma^{n,2} ) n^{2H-1} (
t_{k-1}-t_{k-2} )^{2H-\varepsilon}\\
&&\qquad\sim\biggl( \frac{1}{n}
\biggr)^{1-\varepsilon} \to0\qquad\mbox{as } n \to\infty.
\end{eqnarray*}
This means that we have proven one of the necessary relations:
$n^{2H-1}\sum_{k=2}^{n} J^{n,1}_k\times\break J^{n,2}_k \stackrel{\PP} \to0$ as $n
\to\infty$.

Consider
\begin{eqnarray*}
&&
n^{2H-1}\sum_{k=2}^{n} J^{n,1}_kJ^{n,2}_k\\
&&\qquad=
n^{2H-1} \sum_{k=2}^{n} \int_0^{t_{k-2}} p_k^t (s) \,d W_s
\\
&&\qquad\quad\hspace*{44.2pt}{}\times\int_{t_{k-1}}^{t_{k}} \biggl( \frac{s}{t_{k}} \biggr)^{1/2-H}
(t_{k}-s )^{H-1/2} \,dW_s.
\end{eqnarray*}
As before, it is sufficient to prove that
\begin{eqnarray}
n^{4H-2} \sum_{k=2}^{n} \biggl(\int_0^{t_{k-2}} p_k^t (s)\,
d{W}_s \biggr)^2 \cdot\int_{t_{k-1}}^{t_{k}} \biggl(
\frac{s}{t_{k}} \biggr)^{1-2H} (t_{k}-s )^{2H-1}
\,d[W]_s \stackrel{\PP} \to0 \nonumber\\
\eqntext{\mbox{as }n \to\infty}
\end{eqnarray}
or, equivalently,
%
%e4.8 ###
\begin{equation}\label{eq:1.18.16}
n^{2H-1} \max_{k}\int_{t_{k-1}}^{t_{k}} (t_{k}-s
)^{2H-1} \,d[W]_s \cdot
(\sigma^{n,1} +\sigma^{n,2} ) \stackrel{\PP} \to0.
\end{equation}

Note that by \cite{nvv} and due to the H\"{o}lder properties
of $[W]$,
\[
\int_{t_{k-1}}^{t_{k}} (t_{k}-s )^{2H-1} \,d[W]_s \leq C(\omega
) (t_{k}-
t_{k-1} )^{2H-\varepsilon} \sim\biggl(
\frac{1}{n} \biggr)^{2H-\varepsilon},
\]
whence we obtain \eqref{eq:1.18.16}.

Now, consider $n^{2H-1} \sum J^{n,1}_kJ^{n,4}_k$; other sums can be
estimated similarly. After some transformations, we obtain
\begin{eqnarray*}
&&
n^{4H-2}
\sum_{k=2}^{n} \biggl( \int_0^{t_{k-2}} p_k^t (v) \,d {W}
_u \biggr)^2 \\
&&\qquad\quad\hspace*{10.6pt}{}\times\int_{t_{k-1}}^{t_{k}} s^{1-2H} \biggl(
\int_s^{t_{k}} u^{H-3/2} (u-s)^{H-1/2} \,du \biggr)^2 \,
d[{W}]_s
\\
&&\qquad\leq n^{2H-1}
\max_{k}\int_{t_{k-1}}^{t_{k}} \biggl( \int_s^{t_{k}} u^{H-3/2}
(u-s)^{H-1/2} \,du \biggr)^2 \,d[{W}]_s
(\sigma^{n,1} +\sigma^{n,2} )
\\
&&\qquad
\leq n^{2H-1} \max_{k}\int_{t_{k-1}}^{t_{k}}
\biggl( \int_s^{t_{k}} u^{2H-3} \,du
\int_s^{t_{k}}(u-s)^{2H-1} \,du \biggr)\,
d[{W}]_s \\
&&\qquad\quad{}\times(\sigma^{n,1} +\sigma^{n,2} )\\
&&\qquad
\leq
C n^{2H-1} \max_{k}\int_{t_{k-1}}^{t_{k}} s^{2H-2} (t_{k}-s
)^{2H} \,d[{W}]_s \cdot
(\sigma^{n,1} +\sigma^{n,2} )
\\
&&\qquad
\leq C n \cdot\frac{1}{n}
\max_{k}\int_{t_{k-1}}^{t_{k}} (t_{k}-s )^{2H-1} \,d
[W]_s \cdot(\sigma^{n,1} +\sigma^{n,2})
\\
&&\qquad
\leq C \max_{k} (t_{k}-
t_{k-1} )^{2H-\varepsilon} \cdot(\sigma^{n,1}
+\sigma^{n,2} ) \to0 \qquad\mbox{as } n \to\infty.
\end{eqnarray*}

%s4.4 ###
\subsection{Upper bounds for $[M]$ and $[W]$}

Due to all previous estimates, we can realize our plan and conclude that
\[
t^{2H-1} (t-s)= \lim_{n \to
\infty} n^{2H-1} \sum_{k = n{s/t}+2}^n ( \Delta
X_{t_{k}} )^2 \geq C_1 t^{4H-2} ([M]_t -
[M]_s ),
\]
that is,
\begin{eqnarray*}
[M]_t- [M]_s &\leq& C_2 t^{1-2H}(t-s)= C_2
( t^{-2H}-s t^{1-2H} ) \\
&\leq& C_2 (
t^{-2H} - s^{-2H} )
\end{eqnarray*}
or
\[
\int_s^t u^{1-2H} \,d[W]_u \leq C_2 \int_s^t
u^{1-2H} \,du.
\]

As before, it follows that $[W]_t$ is absolutely
continuous with respect to Lebesgue measure,
%
%e4.9 ###
\begin{equation}\label{eq:1.18.17}
[W]_t = \int_0^t \theta_s \,ds,
\end{equation}
$0 \leq\theta_s \leq C$, where $C$ is some constant and $\theta_s$
is possibly
random. Of course, this is not our final goal, but we can now proceed
with the above estimates for $n^{2H-1} \sum_{k
= n{s/t}+2}^n (J^{n,i} )^2$, $i>1$, and this, together
with condition
(b) [or \eqref{eq:vottak}] and \eqref{eq:half-a-3}, will give us the
possibility to obtain a lower
bound for $[W]_t-[W]_s$, that is, to obtain \eqref{eq:tochtonado}.

%s4.5 ###
\subsection{Lower bound for $[W]_t-[W]_s$}
We can continue estimating from above: for
example, if we take, for simplicity, the sums over $k=2$ up to $k=n$,
then
\begin{eqnarray*}
&&
n^{2H-1} \sum_{k = 2}^n ( J^{n,2}_k )^2 \\
&&\qquad=
\widetilde{\sigma}^{n,1} + \widetilde{\sigma}^{n,2}
\\
&&\qquad:= C
n^{2H-1}\sum_{k = 1}^n \int_{t_{k-2}}^{t_{k-1}}
\biggl(\int_{t_{k-1}}^{t_{k}} \biggl( \frac{s}{u} \biggr)^{1/2-H}
(u-s)^{H-3/2} \,du \biggr)^2\,d[W]_s\\
&&\qquad\quad{} + C n^{2H-1} \sum_{k = 1}^n
\int_{t_{k-2}}^{t_{k-1}} \biggl(\int_{t_{k-2}}^u p_k^t (v) \,dW_v
\biggr) p_k^t(u) \,dW_u
\end{eqnarray*}
and we now need an estimate $ p_k^t (s) \leq( t_{k-1} -s
)^{H-1/2} C.
$

Therefore,
\[
\widetilde{\sigma}^{n,1} \leq C n^{2H-1} \sum_{k = 1}^n
\int_{t_{k-2}}^{t_{k-1}} ( t_{k-1} -s )^{2H-1} \,d
[W]_s.
\]

%Direct estimates now give nothing
We cannot now continue to
estimate the last expression directly
(because of the singularity at the
upper point $t_{k-1}$).
So, we take an indirect route: for some $A>0$,
\begin{eqnarray*}
&&\int_{t_{k-2}}^{t_{k-1}} ( t_{k-1} -s )^{2H-1} \,d[W]_s \\
&&\qquad\leq \int_{t_{k-2}}^{t_{k-1} - {t}/({nA})}
+ \int_{t_{k-1} - {t}/({nA})}^{t_{k-1}}
\\
&&\qquad\leq
\bigl( t_{k-1} - \bigl(t_{k-1} - {t}/({nA}) \bigr) \bigr)^{2H-1}
\cdot\Delta[W]_{t_{k}}
\\
&&\qquad\quad{} + \mbox{[thanks to
\eqref{eq:1.18.17}] } C \int_{t_{k-1} - {t}/({nA})}^{t_{k-1}}
( t_{k-1} -s )^{2H-1} \,ds
\\
&&\qquad\leq\biggl( \frac{t}{nA} \biggr)^{2H-1} \Delta[W
]_{t_{k}} + C \biggl( \frac{t}{nA} \biggr)^{2H}.
\end{eqnarray*}
Taking the
sum, we obtain
\begin{eqnarray*}
\widetilde{\sigma}^{n,1} &\leq& C n^{2H-1} \sum_{k
= 1}^n \biggl( \frac{t}{nA} \biggr)^{2H-1} \Delta[W
]_{t_{k}} + C n^{2H-1} n \biggl(
\frac{t}{nA} \biggr)^{2H}
\\
&\leq& C A^{1-2H}
t^{2H-1} [W]_{t} + C \frac{1}{A^{2H}}
t^{2H}.
\end{eqnarray*}
If we estimate the sum from $k=n\frac{s}{t} + 1$ to
$k=n$, then
\begin{eqnarray*}
\widetilde{\sigma}^{n,1} &\leq& C A^{1-2H}
t^{2H-1} ([W]_{t} - [W
]_{s} ) + C \frac{1}{A^{2H}} t^{2H} \biggl(
1- \frac{s}{t} \biggr)
\\
&=& C A^{1-2H} t^{2H-1} ([
W]_{t} - [W]_{s} ) + C
\frac{1}{A^{2H}} t^{2H-1} (t-s).
\end{eqnarray*}

We now want to prove that $\widetilde{\sigma}^{n,2} \stackrel{\PP}
\to
0$ as $n \to\infty.
$
As usual, it is enough to establish that
\[
n^{4H-2} \sum_{k=1}^{n} \int_{t_{k-2}}^{t_{k-1}}
\biggl(\int_{t_{k-2}}^u p_k^t (v) \,d{W}_v \biggr)^2
( p_k^t(u) )^2 \,d [W]_u
\stackrel{\PP} \to0.
\]
We can now bound $[W]_u$ by $C \,du$,
take the mathematical expectation and note that $ ( p_k^t(u)
)^2 \leq C n^{1-2H}$.
Therefore, it is sufficient to prove that
\[
n^{4H-2} \sum_{k=1}^{n} \int_{t_{k-2}}^{t_{k-1}} \int_{t_{k-2}}^u
(p_k^t (v) )^2 \,d[{W}]_v (
p_k^t(u) )^2 \,du \stackrel{\PP} \to0.
\]
Since $C \,dv$ bounds $d[
{W}]_v$, we have that this value can be bounded by
\begin{eqnarray*}
&&
Cn^{4H-2} \sum_{k=1}^{n} \int_{t_{k-2}}^{t_{k-1}}
\biggl(\int_{t_{k-2}}^u (p_k^t (v) )^2 \,dv \biggr)
( p_k^t(u) )^2 \,du
\\
&&\qquad
\leq C\sum_{k=1}^{n} \int_{t_{k-2}}^{t_{k-1}}
\biggl(\int_{t_{k-2}}^u \,dv \biggr) \,du \leq\frac{1}{n}C \to0 \qquad\mbox{as }
n \to\infty.
\end{eqnarray*}

Finally,
\[
n^{2H-1} \sum_{k = n{s/t}+2}^n (
J^{n,2}_k )^2 \leq C A^{1-2H} t^{2H-1} ([W
]_t - [W]_s )+ C \frac{1}{A^{2H}}
t^{2H-1} (t-s).
\]

Now, proceed with $J^{n,3}_k$:
\begin{eqnarray*}
&&
n^{2H-1} \sum_{k = 1}^n (
J^{n,3}_k )^2 \\
&&\qquad= n^{2H-1} \sum_{k = 1}^n \int_{t_{k-1}}^{t_{k}}
\biggl( \biggl(\frac{s}{t_{k}} \biggr)^{1/2-H} (t_{k}-s
)^{H-1/2} \biggr)^2 \,d[W]_s\\
&&\qquad\quad{}+n^{2H-1}
\sum_{k = 1}^n \int_{t_{k-1}}^{t_{k}} \biggl( \int_{t_{k-1}}^u \biggl(
\frac{s}{t_{k}} \biggr)^{1/2-H} ( t_{k}-s )^{H-1/2} \,dW_s
\biggr)\\
&&\qquad\quad\hspace*{78.1pt}{}\times \biggl(\frac{u}{t_{k}} \biggr)^{1/2-H} (t_{k}-u
)^{H-1/2} \,dW_u.
\end{eqnarray*}
The first term can be estimated as
\begin{eqnarray*}
&&
n^{2H-1} \sum_{k = 1}^n
\int_{t_{k-1}}^{t_{k}} ( t_{k}-s )^{2H-1} \,d[
W]_s \\
&&\qquad\leq C \biggl(\frac{t}{A} \biggr)^{2H-1} ( [
W]_t - [W]_s ) + \frac{C}{A^{2H}}
t^{2H-1} (t-s)
\end{eqnarray*}
as before.

And, with the bound $d [W]_s \leq C \,ds$, the second term
can be estimated as $ n^{4H-2} \sum_{k=1}^{n}
\int_{t_{k-1}}^{t_{k}} \int_{t_{k-1}}^u (
t_{k}-s )^{2H-1} \,ds \cdot( t_{k}-u )^{2H-1}
\,du %$$\leq n^{4\alpha} \sum_{k=1}^{n} (\int_{t_{k-1}}^{t_{k}}
% ( t_{k}-s )^{2\alpha} \,ds )^2 \leq C n^{4\alpha}
%$$ $
\leq C n^{-2} \to0$. Therefore, for
$\sum( J^{n,3}_k )^2$, we have the same estimate as for $\sum
( J^{n,2}_k )^2$. Finally, estimate
\begin{eqnarray*}
&&
n^{2H-1} \sum_{k =
1}^n ( J^{n,4}_k )^2 \\
&&\qquad= C n^{2H-1} \sum_{k=1}^n \biggl(
\int_{t_{k-1}}^{t_{k}} s^{1/2-H} \int_s^{t_{k}}
u^{H-3/2}(u-s)^{H-1/2}\,du \,dW_s \biggr)^2
\\
&&\qquad
=C n^{2H-1}
\sum_{k=1}^n \int_{t_{k-1}}^{t_{k}} s^{1-2H} \biggl(
\int_s^{t_{k}} u^{H-3/2}(u-s)^{H-1/2}\,du \biggr)^2 \,
d[W]_s
\\
&&\qquad\quad{}
+C n^{2H-1} \sum_{k=1}^n \int_{t_{k-1}}^{t_{k}}
\int_{t_{k-1}}^u s^{1/2-H} \int_s^{t_{k}}
v^{H-3/2}(v-s)^{H-1/2}\,dv \,dW_s
\\
&&\qquad\quad\hspace*{109.3pt}{}
\times u^{1/2-H} \int_u^{t_{k}} v^{H-3/2}(v-u)^{H-1/2}\,dv\,
dW_u.
\end{eqnarray*}

The first term can be estimated with the help of \eqref{eq:1.18.17}
as
\begin{eqnarray*}
&&
n^{2H-1} t^{1-2H} \sum_{k=2}^n \int_{t_{k-1}}^{t_{k}} \biggl( \int
_s^{t_{k}} u^{H-3/2}
(u-s)^{H-1/2} \,du \biggr)^2 \,d[W]_s
%$$\leq n^{2\alpha} t^{-2\alpha} \sum_{k=2}^n \int_{t_{k-1}}^{t_{k}}
%$$\leq C n^{2\alpha} t^{-2\alpha} \sum_{k=2}^n \int_{t_{k-1}}^{t_{k}}
%s^{2\alpha-1} ( t_{k}-s
% )^{2\alpha+1} d[ W \rangle_s $$
%$$\leq C
%n^{2\alpha}t^{-2\alpha}\sum_{k=2}^{n} (t_{k-1} )^{2\alpha-1}
% (\frac1n )^{2\al+2}\leq C
%n^{2\alpha-2H}\int_{1/n}^1s^{2\alpha-1}\,ds$$
\\
&&\qquad\leq Cn^{-2H}\to0 \qquad\mbox{as } n\to\infty.
\end{eqnarray*}

If $k=1$, then, for $\frac1p+\frac1q=1$, $p,q>1$,
\begin{eqnarray*}
&&
n^{2H-1}
t^{1-2H} \int_0^{t/n} \biggl( \int_s^{t/n} u^{H-3/2}
(u-s)^{H-1/2}\,du \biggr)^2 \,ds
\\
&&\qquad\leq n^{2H-1} t^{1-2H}
\int_0^{t/n} \biggl( \int_s^{t/n}u^{p(H-3/2)}\,du \biggr)^{2/p}\\
&&\qquad\quad\hspace*{78.1pt}{}\times
\biggl(\int_s^{t/n} (u-s)^{(H-1/2)q} \,du \biggr)^{2/q} \,ds
\\
&&\qquad\leq
n^{2H-1} t^{1-2H} \int_0^{t/n}
s^{ (pH-{3p}/{2}+1 ){2/p}} \biggl(
\frac{t}{n}-s \biggr)^{ ( Hq- {q}/{2}+1 ) {2/q}}
\,ds
\\
&&\qquad
=n^{2H-1} t^{1-2H} \int_0^{t/n}
s^{2H-3+{2/p}} \biggl( \frac{t}{n}-s \biggr)^{2H-1+
{2/q}} \,ds\\
&&\qquad\sim n^{2H-1} t^{1-2H} \biggl(
\frac{t}{n} \biggr)^{4H-1} \to0,
\end{eqnarray*}
that is, the ``main term'' of
$n^{2H-1} \sum_{k=1}^n ( J^{n,4}_k )^2$ tends to 0. For the
remainder term of $n^{2H-1} \sum_{k=1}^n ( J^{n,4}_k )^2$,
it is sufficient to prove that for any $\varepsilon>0$,
\begin{eqnarray}
&&\widetilde{\sigma}^{n,3}:= n^{4H-2}\sum_{k={n\varepsilon}/t}^n
\int_{t_{k-1}}^{t_{k}} \int_{t_{k-1}}^{u} \biggl( s^{1/2-H}
\int_s^{t_{k}} v^{H-3/2} (v-s)^{H-1/2}\,dv \biggr)^2 \,ds
\nonumber\\
&&\qquad\quad\hspace*{101.6pt}\hspace*{-22.21pt}
{}\times u^{1-2H} \biggl(\int_u^{t_{k}} v^{H-3/2}
(v-u)^{H-1/2}\,dv \biggr)^2\,du \to0\nonumber\\
\eqntext{\mbox{as } n\to\infty.}
\end{eqnarray}
However,
\begin{eqnarray*}
\widetilde{\sigma}^{n,3} &\leq& n^{4H-2}
\sum^{n}_{k={n\varepsilon}/t} \int_{t_{k-1}}^{t_{k}}
\int_{t_{k-1}}^{u} \biggl( \int_s^{t_{k}} v^{H-3/2}
(v-s)^{H-1/2}\,dv \biggr)^2 \,ds
\\
&&\hspace*{99.5pt}\hspace*{-23.1pt}{}
\times\biggl( \int_u^{t_{k}}
v^{H-3/2} (v-u)^{H-1/2}\,dv \biggr)^2\,du \\
&\leq&
n^{-6}\sum^{n}_{k={n\varepsilon}/t} (t_{k-1} )^{-4}\sim
n^{-2}\to0 \qquad\mbox{as } n\to\infty.
\end{eqnarray*}
After all estimates, for $s >0$,
\begin{eqnarray*}
&&
\lim_{n \to\infty}n^{2H-1}
\sum_{k=n{s/t}+2}^n ( \Delta X_{t_{k,n}} )^2 \\
&&\qquad\leq C_2
A^{1-2H} t^{2H-1} ([W]_t - [W
]_s ) + C_2 \frac{1}{A^{2H}} t^{2H-1}(t-s).
\end{eqnarray*}
We have the opposite estimate,
\begin{eqnarray*}
C_1 t^{2H-1} (t-s) &\leq& \lim_{n \to\infty}n^{2H-1}
\sum_{k=n{s/t}+2}^n ( \Delta X_{t_{k,n}} )^2
\\
&\leq& C_2 A^{1-2H} t^{2H-1} ([W]_t -
[W]_s )+ C_2 \frac{1}{A^{2H}}
t^{2H-1}(t-s).
\end{eqnarray*}
So, for $A$ sufficiently large, $C_3 := C_1 - C_2
\frac{1}{A^{2H}}
>0$, and we obtain that
\[
C_3 t^{2H-1} (t-s) \leq C_2 A^{1-2H} t^{2H-1}
([W]_t - [W]_s ),
\]
whence $[W]_t - [W]_s \geq
\frac{C_3}{C_2} A^{2H-1} (t-s)$, where the constants do not depend on
$s$ and $t$. Therefore, if we write $[W]_t =
\int_{0}^t p_s \,ds$, then $\varepsilon_1 \leq p_s \leq
\varepsilon_2$, $\varepsilon_i >0$ and $ W_t = \int_0^t
p_s^{1/2}\,dV_s$ for some Wiener process $V$. We can then complete the
proof of the theorem using the same arguments as for $H \in(1/2,1)$.

\section*{Acknowledgments}
The authors are grateful to the anonymous referees for careful and
constructive reading, and for useful suggestions.

% imsref loaded by lrinkeviciute, 2010-08-11 14:53:50
%

%
\printaddresses

\end{document}